\documentclass[12pt,oneside,reqno]{amsart}

\usepackage{txfonts}
\usepackage{bbm}
\usepackage{cases}
\usepackage{amsmath}
\usepackage{graphicx}
\usepackage{mathrsfs}
\usepackage{stmaryrd}
\usepackage{amsfonts}
\usepackage{enumerate,amsmath,amssymb,amsthm}

\numberwithin{equation}{section}

\newcommand{\be}{\begin{eqnarray}}
\newcommand{\ee}{\end{eqnarray}}
\newcommand{\ce}{\begin{eqnarray*}}
\newcommand{\de}{\end{eqnarray*}}
\newtheorem{theorem}{Theorem}[section]
\newtheorem{lemma}[theorem]{Lemma}
\newtheorem{remark}[theorem]{Remark}
\newtheorem{definition}[theorem]{Definition}
\newtheorem{proposition}[theorem]{Proposition}
\newtheorem{Examples}[theorem]{Example}
\newtheorem{corollary}[theorem]{Corollary}

\def\eps{\varepsilon}

\def\p{\partial}

\def\[{{\Big[}}
\def\]{{\Big]}}
\def\<{{\langle}}
\def\>{{\rangle}}
\def\({{\Big(}}
\def\){{\Big)}}

\def\bx{{\mathbf{x}}}

\def\dif{{\mathord{{\rm d}}}}
\def\dis{{\mathord{{\rm \bf d}}}}

\def\no{\nonumber}
\def\={&\!\!=\!\!&}
\def\bt{\begin{theorem}}
\def\et{\end{theorem}}
\def\bl{\begin{lemma}}
\def\el{\end{lemma}}
\def\br{\begin{remark}}
\def\er{\end{remark}}
\def\bd{\begin{definition}}
\def\ed{\end{definition}}
\def\bp{\begin{proposition}}
\def\ep{\end{proposition}}
\def\bc{\begin{corollary}}
\def\ec{\end{corollary}}
\def\bx{\begin{Examples}}
\def\ex{\end{Examples}}

\def\cF{{\mathcal F}}

\def\cJ{{\mathcal J}}

\def\cL{{\mathcal L}}

\def\cT{{\mathcal T}}

\def\mC{{\mathbb C}}

\def\mE{{\mathbb E}}

\def\mH{{\mathbb H}}
\def\mI{{\mathbb I}}

\def\mL{{\mathbb L}}

\def\mN{{\mathbb N}}

\def\mP{{\mathbb P}}

\def\mR{{\mathbb R}}

\def\mV{{\mathbb V}}
\def\mW{{\mathbb W}}

\def\sB{{\mathscr B}}
\def\sC{{\mathscr C}}

\def\sF{{\mathscr F}}

\def\sI{{\mathscr I}}

\def\sL{{\mathscr L}}

\def\sT{{\mathscr T}}

\def\geq{\geqslant}
\def\leq{\leqslant}

\def\div{\mathord{{\rm div}}}

\def\v{{\mathrm v}}
\def\e{{\mathrm e}}

\allowdisplaybreaks

\begin{document}

\title{Fundamental solutions of nonlocal H\"ormander's operators
}

\author{XICHENG ZHANG
\\
\\
{\it Dedicated to the memory of Professor Paul Malliavin}}




\address{
School of Mathematics and Statistics,
Wuhan University, Wuhan, Hubei 430072, P.R.China\\
Email: XichengZhang@gmail.com
 }

\begin{abstract}
Consider the following nonlocal integro-differential operator: for $\alpha\in(0,2)$,
$$
\cL^{(\alpha)}_{\sigma,b} f(x):=\mbox{p.v.} \int_{|z|<\delta}\frac{f(x+\sigma(x)z)-f(x)}{|z|^{d+\alpha}}\dif z+b(x)\cdot\nabla f(x)+\sL f(x),
$$
where $\sigma:\mR^d\to\mR^d\times\mR^d$ and $b:\mR^d\to\mR^d$ are two $C^\infty_b$-functions, $\delta$ is a small positive number, 
p.v. stands for the Cauchy principal value, and $\sL$ is a bounded linear operator in Sobolev spaces.  
Let $B_1(x):=\sigma(x)$ and $B_{j+1}(x):=b(x)\cdot\nabla B_j(x)-\nabla b(x)\cdot B_j(x)$ for $j\in\mN$. 
Under the following uniform H\"ormander's type condition: for some $j_0\in\mN$,
$$
\inf_{x\in\mR^d}\inf_{|u|=1}\sum_{j=1}^{j_0}|u B_j(x)|^2>0,
$$
by using Bismut's approach to the Malliavin calculus with jumps, we prove the existence of fundamental solutions to operator $\cL^{(\alpha)}_{\sigma,b}$. 
In particular, we answer a question proposed by Nualart \cite{Nu1} and Varadhan \cite{Va}.
\end{abstract}

\maketitle \rm

\section{Introduction}

Consider the following nonlocal integro-differential operator: for $\alpha\in(0,2)$,
\begin{align}
\cL^{(\alpha)}_{\sigma,b} f(x):=\mathrm{p.v.}\int_{\mR^d_0}\frac{f(x+\sigma(x)z)-f(x)}{|z|^{d+\alpha}}\dif z+b(x)\cdot\nabla f(x),\label{Op1}
\end{align}
where $\mR^d_0:=\mR^d-\{0\}$ and $\sigma:\mR^d\to\mR^d\times\mR^d$ and $b:\mR^d\to\mR^d$ are two smooth functions and have bounded derivatives of all orders. 
Define $B_1(x):=\sigma(x)$ and $B_{j+1}(x):=b(x)\cdot\nabla B_j(x)-B_j(x)\cdot\nabla b(x)$ for $j\in\mN$.
Recently, in a previous work \cite{Zh4}, we have proved that if for each $x\in\mR^d$, there is a $n(x)\in\mN$ such that
$$
\mathrm{Rank}[B_1(x),\cdots, B_{n(x)}(x)]=d,
$$
then the heat kernel $\rho_t(x,y)$ of operator $\cL^{(\alpha)}_{\sigma,b}$ exists, and as a function of $y$, it is continuous in $L^1(\mR^d)$ with respect to $t,x$. Moreover, 
when $\sigma(x)=\sigma$ is {\it constant}, under the following uniform H\"ormander's type condition: for some $j_0\in\mN$,
$$
\inf_{x\in\mR^d}\inf_{|u|=1}\sum_{j=1}^{j_0}|u B_j(x)|^2>0,
$$
the smoothness of $(t,x,y)\mapsto \rho_t(x,y)$ is also obtained. 
The proofs in \cite{Zh4} are based on the Malliavin calculus to the subordinated Brownian motion (cf. \cite{Ku}), i.e., to consider the following 
stochastic differential equation (abbreviated as SDE):
\begin{align}
\dif X_t(x)=b(X_t(x))\dif t+\sigma(X_{t-}(x))\dif W_{S_t},\ \ X_0(x)=x,\label{SDE0}
\end{align}
where $W_{S_t}$ is an $\alpha$-subordinated Brownian motion. It is well-known that the generator of Markov process $X_t(x)$ is given by $\cL^{(\alpha)}_{\sigma,b}$.
Thus, the main purpose is to study the existence and smoothness of the distribution density $\rho_t(x,y)$ of $X_t(x)$.
If $\sigma(x)$ depends on $x$, since the solution $\{X_t(x),x\in\mR^d\}$ of SDE (\ref{SDE0})
does not form a stochastic diffeomorphism flow in general (cf. \cite{Pr}), it seems impossible to prove the smoothness of $\rho_t(x,y)$ in the framework of \cite{Zh4}.
In this work, we shall study the smoothness of $\rho_t(x,y)$ for {\it non-constant} coefficient $\sigma(x)$ in a different framework. 

As far as we know, Bismut \cite{Bi0} first used Girsanov's transformation to study the smoothness of distribution densities to SDEs with jumps. 
Later, in the monograph \cite{Bi-Gr-Ja}, Bichteler, Gravereaux and Jacod systematically developed the Malliavin calculus with jumps and
studied  the smooth density for SDEs driven by nondegenerate jump noises.
In \cite{Pi}, Picard used difference operator to give another criterion for the smoothness of 
the distribution density of Poisson functionals, and also applied it to SDEs driven by pure jump L\'evy processes. 
By combining the classical Malliavin calculus and Picard's difference operator argument, 
Ishikawa and Kunita in \cite{Is-Ku} obtained a new criterion for the smooth density of Wiener-Poisson functionals (see also \cite{Ku00}).
On the other hand, Cass  in \cite{Ca} established a H\"ormander's type theorem for SDEs with jumps
by proving a Norris' type lemma for discontinuous semimartingales,  but the Brownian diffusion term can not disappear. 
In the pure jump degenerate case, by using a Komatsu-Takeuchi's estimate proven in \cite{Ko-Ta1} for discontinuous semimartingales, 
Takeuchi \cite{Ta} and Kunita \cite{Ku01, Ku0} also obtained similar H\"ormander's theorems. However, their results do not cover operator (\ref{Op1}).
More discussions about their results can be found in \cite{Zh4}.

Let us now consider the following nonlocal integro-differential operator
\begin{align}
\cL_0 f(x):=\mathrm{p.v.}\int_{\Gamma^\delta_0}(f(x+\sigma(x,z))-f(x))\nu(\dif z)+b(x)\cdot\nabla f(x),\label{Op3}
\end{align}
where $\Gamma^\delta_0:=\{0<|z|<\delta\}$, $\sigma(x,z):\mR^d\times\Gamma^\delta_0\to\mR^d$, $b:\mR^d\to\mR^d$, $\nu(\dif z)$ is a L\'evy measure
on $\Gamma^\delta_0$ and satisfies
\begin{align}
\int_{\eps<|z|<\delta}\sigma(x,z)\nu(\dif z)=0,\ \ \forall\eps\in(0,\delta),\ \ x\in\mR^d,\label{14}
\end{align}
and p.v. stands for the Cauchy principal value. Notice that (\ref{14}) is a common assumption in the study of non-local operators.

We make the following assumptions:
\begin{enumerate}[{\bf (H$^\sigma_b$)}]
\item $b$ and $\sigma$ are smooth and for any $m\in\{0\}\cup\mN, j\in\mN$ and some $C_{m}, C_{mj}\geq 1$,
\begin{align}
|\nabla^mb(x)|\leq C_{m},\ \ |\nabla_x^m\nabla^j_z \sigma(x,z)|\leq C_{mj}|z|^{(1-j)\vee 0}.\label{Con1}
\end{align}
\end{enumerate}
\begin{enumerate}[{\bf (H$^\nu$)}]
\item $\nu(\dif z)|_{\Gamma^\delta_0}=\kappa(z)\dif z|_{\Gamma^\delta_0}$, 
$\kappa\in C^\infty(\Gamma^\delta_0;(0,\infty))$ satisfies the following order condition: 
\begin{align}
\lim_{\eps\downarrow 0}\eps^{\alpha-2}\int_{|z|\leq \eps}|z|^2\kappa(z)\dif z=:c_1>0,\label{ET11}
\end{align}
and bounded condition: for $m\in\mN$ and some $C_m\geq 1$,
\begin{align}
|\nabla^m\log\kappa(z)|\leq C_m|z|^{-m}, \ \ z\in \Gamma^{\frac{\delta}{2}}_0.\label{ET01}
\end{align}
\end{enumerate}
\begin{enumerate}[{\bf (UH)}]
\item Let $B_1(x):=\nabla_z\sigma(x,0)$ and define $B_{j+1}(x):=b(x)\cdot\nabla B_j(x)-B_j(x)\cdot\nabla b(x)$ for $j\in\mN$.
Assume that for some $j_0\in\mN$,
\begin{align}
\inf_{x\in\mR^d}\inf_{|u|=1}\sum_{j=1}^{j_0}|uB_j(x)|^2=:c_0>0.\label{UH}
\end{align}
\end{enumerate}
The aim of this paper is to prove the following H\"ormander's type theorem.
\bt\label{Main}
Let $\sL$ be a bounded linear operator in Soblev spaces $\mW^{k,p}(\mR^d)$ for any $p>1$ and $k\in\{0\}\cup\mN$. Under {\bf (H$^\sigma_b$), (H$^\nu$)} and {\bf (UH)}, 
if $\delta\leq\frac{1}{2C_{10}}$, where $C_{10}$ is the same as in (\ref{Con1}), then there exists a measurable function $\rho_t(x,y)$ 
on $(0,1)\times\mR^d\times\mR^d$ called fundamental solution of operator $\cL_0+\sL$ with the properties that
\begin{enumerate}[(i)]
\item For each $t\in(0,1)$ and almost all $y\in\mR^d$, $x\mapsto\rho_t(x,y)$ is smooth, and there is a $\gamma=\gamma(\alpha,j_0,d)$ such that 
for any $p\in(1,\infty)$ and $k\in\{0\}\cup\mN$,
\begin{align}
\|\nabla^k\rho_t(x,\cdot)\|_p\leq C t^{-(k+d)\gamma},\ \ t\in(0,1).\label{NG2}
\end{align}
\item For any $p\in(1,\infty)$ and $\varphi\in L^p(\mR^d)$, 
$\cT_t\varphi(x):=\int_{\mR^d}\varphi(y)\rho_t(x,y)\dif y\in \cap_{k}\mW^{k,p}(\mR^d)$ satisfies
\begin{align}
\p_t \cT_t\varphi(x)=(\cL_0+\sL)\cT_t\varphi(x), \ \ \forall (t,x)\in(0,1)\times\mR^d.\label{NG1}
\end{align}
\end{enumerate} 
\et
\br
The role of operator $\sL$ is usually referred to the large jump part as shown in Corollary \ref{Cor1} below. 
\er

Let us briefly introduce the strategy of proving this theorem. Let $N(\dif t,\dif z)$ be a Poisson random measure with intensity $\dif t\nu(\dif z)$,
and $\tilde N(\dif t,\dif z):=N(\dif t,\dif z)-\dif t\nu(\dif z)$ the compensated Poisson random measure. Consider the following SDE:
$$
X_t(x)=x+\int^t_0b(X_s(x))\dif s+\int^t_0\!\!\!\int_{\Gamma^\delta_0}\sigma(X_{s-},z)\tilde N(\dif s,\dif z).
$$
Under {\bf (H$^\sigma_b$)}, it is well-known that the above SDE has a unique solution $X_t(x)$, which defines a Markov process with generator $\cL_0$.
Let $\cT^0_t f(x):=\mE f(X_t(x))$. Our first aim is to show that under (\ref{Con1})-(\ref{UH}), $X_t(x)$ admits a smooth density, which is a consequence of the
following gradient type estimate: for any $m,k\in\{0\}\cup\mN$, $p\in(1,\infty]$ and $f\in L^p(\mR^d)$,
\begin{align}
\|\nabla^m\cT^0_t\nabla^k f\|_p\leq Ct^{-\gamma_{mk}}\|f\|_p,\label{Short}
\end{align}
where $\nabla$ stands for the gradient operator and $\gamma_{mk}>0$. This will be realized by using Bismut's approach to the Malliavin calculus with jumps.
Of course, the core task is to prove the $L^p$-integrability of the inverse of the reduced Malliavin matrix (see Section 3). 
In order to treat operator $\cL_0+\sL$, letting $\cT_t$ be the corresponding semigroup, by Duhamel's formula, we have
$$
\cT_t\varphi=\cT^0_t\varphi+\int^t_0\cT^0_{t-s}\sL\cT_s\varphi\dif s.
$$
Using the short-time estimate (\ref{Short}) and suitable interpolation techniques, we can prove similar gradient estimates for $\cT_t\varphi$, which
shall produce the desired results by Sobolev's embedding theorem.

As an application of Theorem \ref{Main}, we consider operator $\cL^{(\alpha)}_{\sigma,b}$ in (\ref{Op1}) with $\sigma$ taking the following form:
\begin{align}
\sigma(x)
=\left(
\begin{array}{lcl}
0_{d_1\times d_1},& 0_{d_1\times d_2}\\
0_{d_2\times d_1},&\sigma_0(x)
\end{array}
\right),\label{EQM1}
\end{align}
where $d_1+d_2=d$ and $\sigma_0(x)$ is a $d_2\times d_2$-matrix-valued invertible function.
\bc\label{Cor1}
Assume that $\sigma_0, b$ are smooth and have bounded partial derivatives of all orders, and {\bf (UH)} holds. If $\sigma_0$ satisfies
\begin{align}
\|\sigma^{-1}_0\|_\infty<\infty,\label{EH11}
\end{align}
then the conclusions in Theorem \ref{Main} hold for operator $\cL^{(\alpha)}_{\sigma,b}$.
\ec
\begin{proof}
Let $\chi_\delta:[0,\infty)\to[0,1]$ be a smooth function  with
$$
\chi_\delta(x)=1,\ \ x\in[0,\tfrac{\delta}{2}],\ \ \chi_\delta(x)=0,\ \ x\in[\delta,\infty).
$$
We can write
$$
\cL^{(\alpha)}_{\sigma,b}f(x)=\cL_0f(x)+\sL f(x),
$$
where
$$
\cL_0 f(x):=\mathrm{p.v.}\int_{\Gamma^\delta_0}\frac{f(x+\sigma(x)z)-f(x)}{|z|^{d+\alpha}}\chi_\delta(|z|)\dif z+b(x)\cdot\nabla f(x)
$$
and
\begin{align}
\sL f(x):=\int_{\mR^d_0}\frac{f(x+\sigma(x)z)-f(x)}{|z|^{d+\alpha}}(1-\chi_\delta(|z|))\dif z.\label{Op45}
\end{align}
{\it Claim:} $\sL$ is a bounded linear operator in Sobolev spaces $\mW^{k,p}(\mR^d)$.
\\
{\it Proof of Claim:} Let $z=(z_1,z_2)$ with $z_1\in\mR^{d_1}$ and $z_2\in\mR^{d_2}$. Define
$$
\xi(x,z):=(z_1,\sigma^{-1}_0(x)z_2)\in\mR^d.
$$
Notice that by (\ref{EH11}), there is a positive constant $c_0>0$ such that for all $x,z$,
$$
c_0|z|\leq |\xi(x,z)|\leq c^{-1}_0|z|.
$$
By the change of variables, we have
\begin{align*}
\sL f(x)&=\int_{\mR^d_0}(f(x+(0,z_2))-f(x))\frac{1-\chi_\delta(|\xi(x,z)|)}{|\xi(x,z)|^{d+\alpha}}\det(\sigma^{-1}_0(x))\dif z\\
&=\int_{|z|>\frac{\delta}{2c_0}}(f(x+(0,z_2))-f(x))\frac{1-\chi_\delta(|\xi(x,z)|)}{|\xi(x,z)|^{d+\alpha}}\det(\sigma^{-1}_0(x))\dif z.
\end{align*}
Thus, by Minkovskii's inequality, we have
$$
\|\sL f\|_p\leq \int_{|z|>\frac{\delta}{2c_0}}\|f(\cdot+(0,z_2))-f(\cdot)\|_p\frac{\|\det(\sigma^{-1}_0)\|_\infty}{(c_0|z|)^{d+\alpha}}\dif z\leq C\|f\|_p.
$$
By the chain rule and cumbersome calculations, one sees that
$$
\|\nabla^k\sL f\|_p\leq C\sum_{j=0}^k\|\nabla^j f\|_p.
$$
The proof of claim is complete.
\\
\\
Moreover, if we let $\kappa(z):=\chi_\delta(|z|)|z|^{-d-\alpha}$, then it is easy to check that (\ref{ET01}) is true. Thus, we can use Theorem \ref{Main} to conclude the proof.
\end{proof}
\br
If $\sigma(x)$ is non-degenerate and satisfies (\ref{EH11}), then the law of solutions to SDE (\ref{SDE0}) has a density, 
which is smooth in the first variable. Even in this case, this result seems to be new as all of the well-known results require that
$x\mapsto x+\sigma(x)z$ is invertible (cf. \cite{Pi, Ba-Cl}).
\er
In Corollary \ref{Cor1}, $\sigma$ is required to take a special form (\ref{EQM1}), which has been used to show the boundedness of $\sL$ defined by (\ref{Op45}) in $L^p$-space.
Without assuming (\ref{EQM1}), due to the non-invertibility of $x\mapsto x+\sigma(x)z$, it is not any more true that operator $\sL$ in (\ref{Op45}) is bounded in $L^p$-space.
Consider the following operator:
\begin{align*}
\cL^\nu_{\sigma,b} f(x):=\mathrm{p.v.}\int_{\mR^d_0}(f(x+\sigma(x,z))-f(x))\nu(\dif z)+b(x)\cdot\nabla f(x),
\end{align*}
where $\sigma, b$ and $\nu$ are as above. Let $\cT_t$ be the corresponding semigroup associated to $\cL^\nu_{\sigma,b}$.
Instead of working in the Sobolev space, if we consider the H\"older space, then we have
\bt\label{Main1}
Under {\bf (H$^\sigma_b$), (H$^\nu$)} and {\bf (UH)}, if $\int_{|z|\geq 1}|z|^q\nu(\dif z)<\infty$ for some $q>0$, then there exists a probability density function 
$\rho_t(x,y)$ such that for any $\varphi\in L^\infty(\mR^d)$, $\cT_t\varphi(x)=\int_{\mR^d}\varphi(y)\rho_t(x,y)\dif y$ belongs to H\"older space $\mC^{q+\eps}$, where
$\eps\in(0,1)$ only depends on $\alpha, j_0,d$ with $\alpha$ from (\ref{ET11}) and $j_0$ from (\ref{UH}). Moreover, if $\alpha<q+\eps$, then 
$\p_t \cT_t\varphi(x)=\cL^\nu_{\sigma,b}\cT_t\varphi(x)$ for all $(t,x)\in(0,1)\times\mR^d$.
\et
\br
Consider operator $\cL^{(\alpha)}_{\sigma,b}$ in (\ref{Op1}). Assume that $\sigma$ and $b$  have bounded derivatives of all orders and satisfy {\bf (UH)}.
Let $X_t(x)$ be the unique solution of SDE (\ref{SDE0}). For any $\varphi\in L^\infty(\mR^d)$, 
by the above theorem, $\cT_t\varphi(x)=\mE\varphi(X_t(x))\in \mC^{\alpha+\eps}$ is a classical solution of equation
$$
\p_t u(t,x)=\cL^{(\alpha)}_{\sigma,b} u(t,x).
$$
Here we do not assume that $\sigma$ has the form (\ref{EQM1}). The price we have to pay is that the regularity of $\cT_t\varphi$ depends on the moment of the L\'evy measure.
\er
This paper is organised as follows: In Section 2, we first recall the Bismut's approach to the Malliavin calculus with jumps, 
and an inequality for discontinuous semimaringales proven in \cite{Zh4} which is originally due to Komatsu-Takeuchi \cite{Ko-Ta1}. 
Moreover, we also prove an estimate for exponential Poisson random integrals. In Section 3, we prove a quantity estimate for the Laplace transform of a reduced Malliavin matrix, 
which is the key step in our proofs and can be read independently. We believe that it can be used to other frameworks such as Picard \cite{Pi}
or Ishikawa and Kunita \cite{Is-Ku}. 
In Section 4, we prove the existence of smooth densities for SDEs without big jumps, which corresponds to the operator $\cL_0$ in (\ref{Op3}). 
In Section 5, we treat big jump part and prove our main Theorems \ref{Main} and \ref{Main1} by interpolation and bootstrap arguments.  
Finally, in Appendix, two technical lemmas are proven.

Before concluding this introduction, we collect some notations and make some conventions.
\begin{itemize}
\item Write $\mN_0:=\mN\cup\{0\}$, $\mR^d_0:=\mR^d-\{0\}$.
\item  $\nabla:=(\p_1,\cdots,\p_d)$ denotes the gradient operator.
\item For a c\`adl\`ag function $f:\mR_+\to\mR^d$, $\Delta f_s:=f_s-f_{s-}$.
\item The inner product in Euclidean space is denoted by $\<x,y\>$ or $x\cdot y$.
\item For $p\in[1,\infty]$, $(L^p(\mR^d),\|\cdot\|_p)$ is the $L^p$-space with respect to the Lebesgue measure.
\item $\mW^{k,p}$: Sobolev space; $\mH^{\beta,p}$: Bessel potential space; $\mH^{\beta,\infty}=\mC^\beta$: H\"older space.
\item For a smooth function $f:\mR^d\to\mR^d$, $(\nabla f)_{ij}:=(\p_j f^i)$ denotes the Jacobian matrix of $f$.
\item $C_0(\mR^d)$: The space of all continuous functions with values vanishing at infinity.
\item $C^k_b(\mR^d)$: The space of all bounded continuous functions with bounded continuous partial derivatives up to $k$-order.
Here $k$ can be infinity.
\item The letters $c$ and $C$ with or without indices will denote unimportant constants, whose values may change in different places.
\end{itemize}

\section{Preliminaries}

\subsection{Bismut's approach to the Malliavin calculus with jumps}
In this section, we first recall some basic facts about Bismut's approach to the Malliavin calculus with jumps (cf. \cite{Bi0, Bi-Gr-Ja} and \cite[Section 2]{So-Zh}).
Let $\Gamma\subset\mR^d$ be an open set containing the original point. Let us define
\begin{align}
\Gamma_0:=\Gamma\setminus\{0\},\ \ \varrho(z):=1\vee\dis(z,\Gamma^c_0)^{-1},\label{Rho}
\end{align}
where $\dis(z,\Gamma^c_0)$ is the distance of $z$ to the complement of $\Gamma_0$. Notice that $\varrho(z)=\frac{1}{|z|}$ near $0$.

Let $\Omega$ be the canonical space of all integer-valued measure $\omega$ on $[0,1]\times\Gamma_0$ with $\mu(A)<+\infty$ for any compact set $A\subset[0,1]\times\Gamma_0$.
Define the canonical process on $\Omega$ as follows:
$$
N(\omega; \dif t,\dif z):=\omega(\dif t,\dif z).
$$
Let $(\sF_t)_{t\in[0,1]}$ be the smallest right-continuous filtration on $\Omega$ such that $N$ is optional.
In the following, we write $\sF:=\sF_1$, and endow $(\Omega,\sF)$ with the unique probability measure $\mP$ such that
$N$ is a Poisson random measure with intensity $\dif t\nu(\dif z)$, where $\nu(\dif z)=\kappa(z)\dif z$ with
\begin{align}
\kappa\in C^1(\Gamma_0;(0,\infty)),\  \int_{\Gamma_0}(1\wedge|z|^2)\kappa(z)\dif z<+\infty,\ \ |\nabla\log\kappa(z)|\leq C\varrho(z),\label{ET1}
\end{align}
where $\varrho(z)$ is defined by (\ref{Rho}). In the following we write
$$
\tilde N(\dif t,\dif z):=N(\dif t,\dif z)-\dif t\nu(\dif z).
$$
Let $p\geq 1$ and $k\in\mN$. We introduce the following spaces for later use.
\begin{itemize}
\item $\mL^1_p$: The space of all predictable processes: $\xi:\Omega\times[0,1]\times\Gamma_0\to\mR^k$ with finite norm:
$$
\|\xi\|_{\mL^1_p}:=\left[\mE\left(\int^1_0\!\!\!\int_{\Gamma_0}|\xi(s,z)|\nu(\dif z)\dif s\right)^p\right]^{\frac{1}{p}}
+\left[\mE\int^1_0\!\!\!\int_{\Gamma_0}|\xi(s,z)|^p\nu(\dif z)\dif s\right]^{\frac{1}{p}}<\infty.
$$
\item $\mL^2_p$: The space of all predictable processes: $\xi:\Omega\times[0,1]\times\Gamma_0\to\mR^k$ with finite norm:
$$
\|\xi\|_{\mL^2_p}:=\left[\mE\left(\int^1_0\!\!\!\int_{\Gamma_0}|\xi(s,z)|^2\nu(\dif z)\dif s\right)^{\frac{p}{2}}\right]^{\frac{1}{p}}
+\left[\mE\int^1_0\!\!\!\int_{\Gamma_0}|\xi(s,z)|^p\nu(\dif z)\dif s\right]^{\frac{1}{p}}<\infty.
$$
\item $\mV_p$:  The space of all predictable processes $\v: \Omega\times[0,1]\times\Gamma_0\to\mR^d$ with finite norm:
$$
\|\v\|_{\mV_p}:=\|\nabla_z\v\|_{\mL^1_p}+\|\v\varrho\|_{\mL^1_p}<\infty,
$$
where $\varrho(z)$ is defined by (\ref{Rho}). Below we shall write
$$
\mV_{\infty-}:=\cap_{p\geq 1}\mV_p.
$$
\item $\mV_0$:  The space of all predictable processes $\v: \Omega\times[0,1]\times\Gamma_0\to\mR^d$
with the following properties: (i) $\v$ and $\nabla_z \v$ are bounded;
(ii) there exists a compact subset $U\subset \Gamma_0$ such that
$$
\v(t,z)=0,\ \ \forall z\notin U.
$$
Moreover, $\mV_0$  is dense in $\mV_p$ for all $p\geq 1$ (cf. \cite[Lemma 2.1]{So-Zh}).
\end{itemize}
Let $C_{\mathrm{p}}^\infty(\mR^m)$ be the class of all smooth functions on $\mR^m$ which together with all the derivatives have at most polynomial growth.
Let $\cF C^\infty_{\mathrm{p}}$ be the class of all Poisson functionals on $\Omega$ with the following form:
$$
F(\omega)=f\big(\omega(g_1),\cdots, \omega(g_m)\big),
$$
where $f\in C_{\mathrm{p}}^\infty(\mR^m)$ and $g_1,\cdots, g_m\in\mV_0$ are non-random, and
$$
\omega(g_j):=\int^1_0\!\!\!\int_{\Gamma_0}g_j(s,z)\omega(\dif s,\dif z).
$$
Notice that
$$
\cF C^\infty_{\mathrm{p}}\subset \cap_{p\geq 1}L^p(\Omega,\sF,\mP).
$$
For $\v\in \mV_{\infty-}$ and $F\in\cF C^\infty_{\mathrm{p}}$, let us define
\begin{align*}
D_\v F=:\sum_{j=1}^{m}(\p_j f)(\cdot)\int^1_0\!\!\!\int_{\Gamma_0}\v(s,z)\cdot \nabla_z g_j(s,z)\omega(\dif s,\dif z),
\end{align*}
where ``$(\cdot)$'' stands for ($\omega(g_1),\cdots, \omega(g_m)$).

We have the following integration by parts formula (cf. \cite[Theorem 2.9]{So-Zh}).
\bt\label{Th21}
Let $\v\in\mV_{\infty-}$ and $p>1$.
The linear operator $(D_\v, \cF C^\infty_{\mathrm{p}})$ is closable in $L^p$. The closure is denoted by
$\mW^{1,p}_\v(\Omega)$, which is a Banach space with respect to the norm:
$$
\|F\|_{\v; 1,p}:=\|F\|_{L^p}+\|D_\v F\|_{L^p}.
$$
Moreover, for any $F\in\mW^{1,p}_\v(\Omega)$, we have
\begin{align}
\mE(D_\v F)=-\mE(F \div(\v)),\label{ER88}
\end{align}
where $\div(\v)$ is defined by
\begin{align}
\div(\v):=\int^1_0\!\!\!\int_{\Gamma_0}\frac{\div(\kappa\v)(s,z)}{\kappa(z)}\tilde N(\dif s,\dif z).\label{ER1}
\end{align}
\et

Below, we shall write
$$
\mW_{\v}^{1,\infty-}(\Omega):=\cap_{p>1}\mW_{\v}^{1,p}(\Omega).
$$
The following Kusuoka and Stroock's formula is proven in \cite[Proposition 2.11]{So-Zh}.
\bp\label{Pr1}
Fix $\v\in\mV_{\infty-}$.  Let $\eta(\omega,s,z):\Omega\times[0,1]\times\Gamma_0\to\mR$ be a measurable map
and satisfy that for each $(\omega, s,z)\in\Omega\times[0,1]\times\Gamma_0$, $\eta(\cdot, s,z)\in\mW^{1,\infty-}_\v(\Omega),\ \ \eta(\omega, s,\cdot)\in C^1(\Gamma_0)$,
and
\begin{align}
\mbox{$s\mapsto \eta(s,z), D_\v\eta(s,z),\nabla_z\eta(s,z)$ are left-continuous and $\sF_s$-adapted},
\end{align}
and  for any $p\geq 1$,
\begin{align}
\mE\left[\sup_{s\in[0,1]}\sup_{z\in\Gamma_0}\left(\frac{|\eta(s,z)|^p+|D_\v\eta(s,z)|^p}{(1\wedge|z|)^p}+|\nabla_z\eta(s,z)|^p\right)\right]<+\infty.\label{ET7}
\end{align}
Then $\sI(\eta):=\int^1_0\!\int_{\Gamma_0}\eta(s,z)\tilde N(\dif s,\dif z)\in\mW^{1,\infty-}_\v(\Omega)$ and
\begin{align}
D_\v \sI(\eta)=\int^1_0\!\!\!\int_{\Gamma_0}D_\v\eta(s,z)\tilde N(\dif s,\dif z)+\int^1_0\!\!\!\int_{\Gamma_0}\<\nabla_z\eta(s,z),\v(s,z)\>N(\dif s,\dif z).\label{For2}
\end{align}
\ep

We also need the following Burkholder's inequalities (cf. \cite[Lemma 2.3]{So-Zh}).
\bl\label{Le2}
\begin{enumerate}[(i)]
\item For any $p>1$, there is a constant $C_p>0$ such that for any $\xi\in \mL^1_p$,
\begin{align}
\mE\left(\sup_{t\in[0,1]}\left|\int^t_0\!\!\!\int_{\Gamma_0}\xi(s,z)N(\dif s,\dif z)\right|^p\right)\leq C_p\|\xi\|^p_{\mL^1_p}.\label{BT1}
\end{align}
\item For any $p\geq 2$, there is a constant $C_p>0$ such that for any $\xi\in \mL^2_p$,
\begin{align}
\mE\left(\sup_{t\in[0,1]}\left|\int^t_0\!\!\!\int_{\Gamma_0}\xi(s,z)\tilde N(\dif s,\dif z)\right|^p\right)\leq C_p\|\xi\|^p_{\mL^2_p}.\label{BT2}
\end{align}
\end{enumerate}
\el

\subsection{Two Lemmas}

We first recall the following important Komatsu-Takeuchi's type estimate proven in \cite[Theorem 4.2]{Zh4}, which will be used in Section 3.
\bl\label{Th42}
Let $(f_t)_{t\geq 0}$ and $(f^0_t)_{t\geq 0}$ be two $m$-dimensional semimartingales given by
\begin{align*}
f_t&=f_0+\int^{t\wedge\tau}_0f^0_s\dif s+\int^t_0\!\!\!\int_{|z|\leq\delta}g_{s-}(z)\tilde N(\dif s,\dif z),\\
f^0_t&=f^0_0+\int^t_0f^{00}_s\dif s+\int^t_0\!\!\!\int_{|z|\leq\delta}g^0_{s-}(z)\tilde N(\dif s,\dif z),
\end{align*}
where $\delta\in(0,1]$, $\tau$ is a stopping time and $f_t, f^0_t, f^{00}_t$ and $g_t(z), g^0_t(z)$ are c\`adl\`ag $\sF_t$-adapted processes. Assume that for some $\kappa\geq 1$,
$$
|f_t|^2\vee|f^0_t|^2\vee\sup_{z}\frac{|g_t(z)|^2\vee |g^0_t(z)|^2}{1\wedge|z|^2}\leq\kappa,\ \ a.s.
$$
Then for any $\eps, T\in(0,1]$, there exists a positive random variable $\zeta$ with $\mE \zeta\leq 1$ such that
\begin{align}
c_0\int^{T\wedge\tau}_0|f^0_t|^2\dif t\leq(\delta^{-\frac{3}{2}}+\eps^{-\frac{3}{2}})\int^T_0|f_t|^2\dif t
+\kappa\delta^{\frac{1}{2}}\log \zeta+\kappa(\eps\delta^{-\frac{1}{2}}+\eps^{\frac{1}{2}}+T\delta^{\frac{1}{2}}),\label{Eso2}
\end{align}
where $c_0\in(0,1)$ only depends on $\int_{|z|\leq 1}|z|^2\nu(\dif z)$.
\el
The following result will be used in Section 4.
\bl\label{Le25}
Let $g_s(z),\eta_s$ be two left continuous $\sF_s$-adapted processes satisfying
\begin{align}
0\leq g_s(z)\leq\eta_s,\ |g_s(z)-g_s(0)|\leq\eta_s |z|,\ \forall |z|\leq 1,\label{UY22}
\end{align}
and for any $p\geq 2$,
$$
\mE\left(\sup_{s\in[0,1]}|\eta_s|^p\right)<+\infty.
$$
If for some $\alpha\in(0,2)$,
\begin{align}
\lim_{\eps\to 0}\eps^{\alpha-2}\int_{|z|\leq\eps}|z|^2\nu(\dif z)=:c_1>0,\label{UY21}
\end{align}
then for any $\delta\in(0,1)$, there exist constants $c_2,\theta\in(0,1), C_2\geq 1$ such that for all $\lambda,p\geq 1$ and $t\in(0,1)$,
\begin{align}
\mE\exp\left\{-\lambda\int^t_0\!\!\!\int_{\mR^d_0}g_s(z)\zeta(z)N(\dif s,\dif z)\right\}
\leq C_2\left(\mE\exp\left\{-c_2\lambda^\theta\int^t_0g_s(0)\dif s\right\}\right)^{\frac{1}{2}}+C_p\lambda^{-p},\label{TR7}
\end{align}
where $\zeta(z)=\zeta_\delta(z)$ is a nonnegative smooth function with
$$
\zeta_\delta(z)=|z|^3,\ \ |z|\leq\delta/4,\ \ \zeta_\delta(z)=0,\ \ |z|>\delta/2.
$$
\el
\begin{proof}
For $\lambda\geq 1$ and $\beta>0$, define a stopping time 
$$
\tau:=\inf\{s>0: \eta_s\geq\lambda^\beta\}\wedge 1.
$$
Set
$$
h^\lambda_t:=\int_{\mR^d_0}(1-\e^{-\lambda g_t(z)\zeta(z)})\nu(\dif z)
$$
and
\begin{align}
M^\lambda_t:=-\lambda\int^{t\wedge\tau}_0\!\!\!\int_{\mR^d_0}g_s(z)\zeta(z)N(\dif s,\dif z)+\int^{t\wedge\tau}_0h^\lambda_s\dif s.\label{EY1}
\end{align}
By It\^o's formula, we have
$$
\e^{M^\lambda_t}=1+\int^{t\wedge\tau}_0\!\!\!\int_{\mR^d_0}\e^{M^\lambda_{s-}}(\e^{-\lambda g_s(z)\zeta(z)}-1)\tilde N(\dif s,\dif z).
$$
Since for any $x\geq 0$,
$$
1-\mathrm{e}^{-x}\leq 1\wedge x,
$$
by (\ref{UY22}) and definition of $\tau$, we have
\begin{align*}
M^\lambda_t\leq \int^{t\wedge\tau}_0h^\lambda_s\dif s
&\leq \int^{t\wedge\tau}_0\!\!\!\int_{\mR^d_0}(1\wedge(\lambda g_s(z)\zeta(z)))\nu(\dif z)\dif s\\
&\leq\int_{\mR^d_0}(1\wedge(\lambda^{1+\beta}\zeta(z)))\nu(\dif z)<\infty.
\end{align*}
Hence, $\mE\e^{M^\lambda_t}=1$ and by (\ref{EY1}) and H\"older's inequality,
\begin{align}
\mE\exp\left\{-\frac{\lambda}{2}\int^{t\wedge\tau}_0\!\!\!\int_{\mR^d_0}g_s(z)\zeta(z)N(\dif s,\dif z)\right\}
\leq\left(\mE\exp\left\{-\int^{t\wedge\tau}_0h^\lambda_s\dif s\right\}\right)^{\frac{1}{2}}.\label{TR2}
\end{align}
Since $1_{s<\tau}g_s(z)\leq\lambda^\beta$ by (\ref{UY22})  and definition of $\tau$, and for any  $x\leq 1$,
$$
1-\mathrm{e}^{-x}\geq \tfrac{x}{\e},
$$
for any $q\geq 1+\beta$, there exists  $c\in(0,1)$ small enough such that for all $\lambda\geq 1$ and $s<\tau$,
\begin{align}
h^{\lambda}_s&\geq\int_{|z|^3\leq c\lambda^{-q}}(1-\e^{-\lambda g_s(z)|z|^3})\nu(\dif z)
\geq\frac{\lambda}{\e}\int_{|z|^3\leq c\lambda^{-q}}g_s(z)|z|^3\nu(\dif z)\no\\
&=\frac{\lambda g_s(0)}{\e}\int_{|z|^3\leq c\lambda^{-q}}|z|^3\nu(\dif z)+\frac{\lambda}{\e}\int_{|z|^3\leq c\lambda^{-q}}(g_s(z)-g_s(0))|z|^3\nu(\dif z).\label{TR1}
\end{align}
Notice that by (\ref{UY21}),  for any $p\geq 2$, there exist constants $c_0,C_0>0$ such that for all $\eps\in(0,1)$ (cf. \cite[Lemma 5.2]{So-Zh}),
\begin{align}
c_0\eps^{p-\alpha}\leq\int_{|z|\leq\eps}|z|^p\nu(\dif z)\leq C_0\eps^{p-\alpha}.\label{UY2}
\end{align}
If we choose
$$
\beta\in (0,\tfrac{\alpha\wedge 1}{3-\alpha}),\ \ 
q=\left\{
\begin{aligned}
&1+\beta, &\alpha\in(0,1], \\
&\tfrac{3(1+\beta)}{4-\alpha},&\alpha\in(1,2),
\end{aligned}
\right.
$$
then by (\ref{TR1}), (\ref{UY22}) and (\ref{UY2}), we further have for all $\lambda\geq 1$ and $s<\tau$,
\begin{align}
h^{\lambda}_s\geq c_2g_s(0)\lambda^{1-\frac{(3-\alpha)q}{3}}-C_1\lambda^{1+\beta-\frac{(4-\alpha)q}{3}}
\geq c_2g_s(0)\lambda^{1-\frac{(3-\alpha)q}{3}}-C_1.\label{TR3}
\end{align}
On the other hand, by Chebyshev's inequality, we have for any $p\geq 2$,
$$
\mP(\tau\leq t)=\mP\left(\sup_{s\in[0,t]}\eta_s>\lambda^\beta\right)\leq\lambda^{-\beta p}\mE\left(\sup_{s\in[0,t]}|\eta_s|^p\right),
$$
which together with (\ref{TR2}) and (\ref{TR3}) yields the desired estimate (\ref{TR7}).
\end{proof}

\section{Estimate of Laplace transform of reduced Malliavin matrix}
This section is devised to be independent of the settings in Subsection 2.1 so that it can be used to other framework such as Picard \cite{Pi}.
Let $L_t$ be a $d$-dimensional L\'evy process with L\'evy measure $\nu$. We assume that the L\'evy measure $\nu$ 
satisfies the following conditions: for some $\alpha\in(0,2)$,
\begin{align}
\int_{|z|<\delta}|z|^2\nu(\dif z)\leq C\delta^{2-\alpha},\ \ \forall \delta\in(0,1), \ \int_{|z|\geq 1}|z|^m\nu(\dif z)<\infty,\ \forall m\in\mN.\label{EUT0}
\end{align}
Let $N(\dif t,\dif z)$ be the Poisson random measure associated with $L_t$, i.e.,
$$
N((0,t]\times E)=\sum_{s\leq t}1_E(\Delta L_s),\ \ E\in\sB(\mR^d_0).
$$
Let $\tilde N(\dif t,\dif z):=N(\dif t,\dif z)-\dif t\nu(\dif z)$ be the compensated Poisson random measure. Consider the following SDE:
\begin{align}
X_t(x)=x+\int^t_0b(X_s(x))\dif s+\int^t_0\!\!\!\int_{\mR^d_0}\sigma(X_{s-}(x),z)\tilde N(\dif s,\dif z),\label{SDE1}
\end{align}
where $b:\mR^d\to\mR^d$ and $\sigma:\mR^d\times\mR^d\to\mR^d$ are two functions satisfying that for any $m\in\mN_0$ and $j=0,1$,
\begin{align}
|\nabla^m b(x)|\leq C,\ |\nabla^m_x\nabla^j_z\sigma(x,z)|\leq C|z|^{1-j}\label{33}
\end{align} 
and
\begin{align}
\int_{r<|z|<R}\sigma(x,z)\nu(\dif z)=0,\ \ 0<r<R<\infty.\label{34}
\end{align} 
Under (\ref{33}), it is well-known that SDE (\ref{SDE1}) has a unique solution denoted by $X_t:=X_t(x)$, 
which defines a $C^\infty$-stochastic flow (cf. \cite{Fu-Ku} and \cite{Pr}). 
Let $J_t:=J_t(x):=\nabla X_t(x)$ be the Jacobian matrix of $X_t(x)$, which solves the following linear matrix-valued SDE:
\begin{align}
J_t=\mI+\int^t_0\nabla b(X_s)J_s\dif s+\int^t_0\!\!\!\int_{\mR^d_0}\nabla_x \sigma(X_{s-},z)J_{s-}\tilde N(\dif s,\dif z).\label{SDE2}
\end{align}
If we further assume 
\begin{align}
\inf_{x\in\mR^d}\inf_{z\in\mR^d}\det (\mI+\nabla_x\sigma(x,z))>0,\label{Con3}
\end{align}
then the matrix $J_t(x)$ is invertible (cf. \cite{Fu-Ku}).
Let $K_t=K_t(x)$ be the inverse matrix of $J_t(x)$. By It\^o's formula, it is easy to see that $K_t$ solves the following linear matrix-valued SDE (cf. \cite[Lemma 3.2]{Zh4}):
\begin{align}
K_t&=\mI-\int^t_0 K_s\nabla b(X_s)\dif s+\int^t_0\!\!\!\int_{\mR^d_0}K_{s-}Q(X_{s-},z)\tilde N(\dif s,\dif z)\no\\
&\qquad-\int^t_0\!\!\!\int_{\mR^d_0}K_{s-}Q(X_{s-},z)\nabla_x\sigma(X_{s-},z)\nu(\dif z)\dif s,\label{SDE3}
\end{align}
where
\begin{align}
Q(x,z):=(\mI+\nabla_x\sigma(x,z))^{-1}-\mI.
\end{align}

First of all, we have the following easy estimate. Since the proof is standard by Gronwall's inequality and Burkholder's inequality, we omit the details.
\bl
Under (\ref{33}) and (\ref{Con3}), we have for any $p\geq 1$,
\begin{align}
\sup_{x\in\mR^d}\mE\left(\sup_{t\in[0,1]}(|J_t(x)|^p+|K_t(x)|^p)\right)<+\infty.\label{EUT1}
\end{align}
\el
We now prove the following crucial lemma.
\bl\label{Le9}
Let $V:\mR^d\to\mR^d\times\mR^d$ be a bounded smooth function with bounded derivatives of all orders.   
Under (\ref{33}), (\ref{34})  and (\ref{Con3}), there exist $\beta_1,\beta_3\in(0,1),\beta_2\geq 1$ only depending on $\alpha$ and
constants $C_1\geq 1$ and $c_1\in(0,1)$ only depending on $b,V$ and $\alpha,\beta_i,\nu$ such that for all $\delta, t\in(0,1)$ and $p\geq 1$,
\begin{align}
&\sup_{|u|=1}\mP\left(\int^t_0|uK_s[b,V](X_s)|^2\dif s\geq t\delta^{\beta_1},\int^t_0|uK_sV(X_s)|^2\dif s\leq t\delta^{\beta_2}\right)
\leq C_1\e^{-c_1t\delta^{-\beta_3}}+C_p\delta^p,\label{EE4}
\end{align}
where $[b,V]:=b\cdot\nabla V-V\cdot\nabla b$.
\el
\begin{proof}
We divide the proof into four steps.
\\
{\bf (1)} Fixing $\delta\in(0,1)$, we decompose the L\'evy process as the small and large jump parts,  i.e., $L_t=L^\delta_t+\hat L^\delta_t$, where
$$
L^\delta_t:=\int_{|z|\leq\delta}z\tilde N((0,t],\dif z), \ \hat L^\delta_t:=\int_{|z|>\delta}zN((0,t],\dif z).
$$
Clearly,
$$
\mbox{$L^\delta_t$ and $\hat L^\delta_t$ are independent}.
$$
Let us fix a path $\hbar$ with finitely many jumps on any finite time interval. Let $X^\delta_t(x;\hbar)$ solve the following SDE:
\begin{align}
X^\delta_t(x;\hbar)&=x+\int^t_0b(X^\delta_s(x;\hbar))\dif s+\sum_{s\leq t}\sigma(X^\delta_{s-}(x;\hbar),\Delta \hbar_s)\no\\
&\quad +\int^t_0\!\!\!\int_{|z|\leq\delta}\sigma(X^\delta_{s-}(x;\hbar),z)\tilde N(\dif s,\dif z).\label{EK3}
\end{align}
Let $K^\delta_t(x;\hbar):=[\nabla X^\delta_t(x;\hbar)]^{-1}$. Clearly, by (\ref{34}) we have
\begin{align}
X_t(x)=X^\delta_t(x;\hbar)|_{\hbar=\hat L^\delta},\ \ K_t(x)=K^\delta_t(x;\hbar)|_{\hbar=\hat L^\delta}.\label{EH2}
\end{align}
Moreover, $K^\delta_t:=K^\delta_t(x;0)$  solves the following equation
\begin{align}
K^\delta_t&=\mI-\int^t_0 K^\delta_s\nabla b(X^\delta_s)\dif s+\int^t_0\!\!\!\int_{|z|\leq\delta}K^\delta_{s-}Q(X^\delta_{s-},z)\tilde N(\dif s,\dif z)\no\\
&\quad-\int^t_0\!\!\!\int_{|z|\leq\delta}K^\delta_{s-}Q(X^\delta_{s-},z)\nabla_x \sigma(X^\delta_{s-},z)\nu(\dif z)\dif s.\label{EK33}
\end{align}
{\bf (2)} Define functions:
\begin{align}
H_V(x,z)&:=V(x+\sigma(x,z))-V(x)+Q(x,z)V(x+\sigma(x,z)),\label{Re12}\\
G_V(x,z)&:=H_V(x,z)+\nabla_x\sigma(x,z)\cdot V(x)-\sigma(x,z)\cdot\nabla V(x)\label{Re22}
\end{align}
and
\begin{align*}
V_0(x)&:=[b,V](x)+\int_{|z|\leq\delta}G_V(x,z)\nu(\dif z),\\
V_1(x)&:=[b,V_0](x)+\int_{|z|\leq\delta}G_{V_0}(x,z)\nu(\dif z).
\end{align*}
It is easy to see by (\ref{33}) that
\begin{align}
|H_V(x,z)|\leq C(1\wedge |z|),\ \ |G_V(x,z)|\leq C(1\wedge |z|^2).\label{EH1}
\end{align}
For a row vector $u\in\mR^d$, we introduce the processes:
$$
f_t:=uK^\delta_t V(X^\delta_t),\  f^0_t:=uK^\delta_t V_0(X^\delta_t),\  f^{00}_t:=uK^\delta_t V_1(X^\delta_t),
$$
$$
g_t(z):=uK^\delta_tH_V(X^\delta_t,z),\ \ g^0_t(z):=uK^\delta_t H_{V_0}(X^\delta_t,z),
$$
where
$$
X^\delta_t:=X^\delta_t(x;0),\ \ K^\delta_t:=K^\delta_t(x;0).
$$
By equations (\ref{EK33}),  (\ref{EK3}) with $\hbar=0$ and using It\^o's formula, we have
\begin{align*}
f_t&=u V(x)+\int^t_0uK^\delta_s[b, V](X^\delta_s)\dif s+\int^t_0\!\!\!\int_{|z|\leq\delta}g_{s-}(z)\tilde N(\dif s,\dif z)\\
&\qquad+\int^t_0\!\!\!\int_{|z|\leq\delta}uK^\delta_sG_V(X^\delta_s,z)\nu(\dif z)\dif s\\
&=u V(x)+\int^t_0f^0_s\dif s+\int^t_0\!\!\!\int_{|z|\leq\delta}g_{s-}(z)\tilde N(\dif s,\dif z)
\end{align*}
and
$$
f^0_t=uV_0(x)+\int^t_0f^{00}_s\dif s+\int^t_0\!\!\!\int_{|z|\leq\delta}g^0_{s-}(z)\tilde N(\dif s,\dif z).
$$
For $\gamma\in(0,\frac{1}{4})$, define a stopping time
$$
\tau:=\tau_u(x):=\inf\big\{s\geq 0: |uK^\delta_s(x;0)|^2>\delta^{-\gamma}\big\}.
$$
Since $V$ has bounded derivatives of all orders, there exists a constant $\kappa_0\geq 1$ only depending on $b, \sigma$ and $V$ such that
for all $t\in[0,\tau)$ and $z\in\mR^d$,
$$
|f_t|^2,\ |f^0_t|^2,\ |f^{00}_t|^2\leq \kappa_0\delta^{-\gamma},\ \ |g_t(z)|^2,\ |g^0_t(z)|^2\leq \kappa_0\delta^{-\gamma}(1\wedge|z|^2).
$$
If we make the following replacement in Theorem \ref{Th42}:
$$
f_t,  g_t(z), f^0_t, g^0_t(z)\Rightarrow f_{t\wedge\tau}, 1_{t<\tau}g_t(z), f^0_{t\wedge\tau},  1_{t<\tau}g^0_t(z),
$$ 
then by (\ref{Eso2}) with $\eps=\delta^5$ and $\kappa=\kappa_0\delta^{-\gamma}$, we obtain
\begin{align*}
c_0\int^{t\wedge\tau}_0|f^0_s|^2\dif s&\leq(\delta^{-\frac{3}{2}}+\delta^{-\frac{15}{2}})\int^t_0|f_{s\wedge\tau}|^2\dif s
+\kappa_0\delta^{\frac{1}{2}-\gamma}\log \zeta+\kappa_0\delta^{-\gamma}(\delta^{5-\frac{1}{2}}+\delta^{\frac{5}{2}}+t\delta^{\frac{1}{2}})\\
&\leq 2\kappa_0\delta^{-\frac{15}{2}}\int^t_0|f_{s\wedge\tau}|^2\dif s+\kappa_0\delta^{\frac{1}{2}-\gamma}\log \zeta
+2\kappa_0(\delta^{\frac{5}{2}-\gamma}+t\delta^{\frac{1}{2}-\gamma})\ \ a.s.,
\end{align*}
where $c_0\in(0,1)$  only depends on $\int_{|z|\leq 1}|z|^2\nu(\dif z)$. From this, dividing both sides by $2\kappa_0\delta^{\frac{1}{2}-\gamma}$ and
taking exponential, then multiplying $1_{\tau\geq t}$ and taking expectations, we derive that
\begin{align}
&\mE\left(\exp\left\{\frac{c_0\delta^{\gamma-\frac{1}{2}}}{2\kappa_0}\int^t_0|f^0_s|^2\dif s-\delta^{\gamma-8}\int^t_0|f_s|^2\dif s\right\}1_{\tau\geq t}\right)\no\\
&\qquad\leq \mE(1_{\tau\geq t}\zeta)\exp\{\delta^2+t\}\leq\exp\{\delta^2+t\}.\label{EC2}
\end{align}
Recalling the definition of $f^0_t$ and by $|x+y|^2\geq\frac{|x|^2}{2}-|y|^2$, we have for $t<\tau$,
\begin{align*}
|f^0_t|^2&\stackrel{(\ref{EH1})}{\geq}\frac{|uK^\delta_t[b,V](X^\delta_t)|^2}{2}-C|uK^\delta_t|^2\left(\int_{|z|\leq\delta}|z|^2\nu(\dif z)\right)^2\\
&\stackrel{(\ref{EUT0})}{\geq} \frac{|uK^\delta_t[b,V](X^\delta_t)|^2}{2}-C_2\delta^{4-2\alpha-\gamma}.
\end{align*}
Thus, by (\ref{EC2}) there exist $c_1\in(0,1)$ and $C_3\geq 1$ independent of the starting point $x$ such that for all $\delta,t\in(0,1)$,
\begin{align}
&\sup_{u\in\mR^d}\mE\left(\exp\left\{c_1\delta^{\gamma-\frac{1}{2}}\int^t_0|uK^\delta_s[b,V](X^\delta_s)|^2\dif s
-\delta^{\gamma-8}\int^t_0|uK^\delta_sV(X^\delta_s)|^2\dif s\right\}1_{\tau\geq t}\right)\no\\
&\qquad\qquad\qquad\leq\exp\Big\{\delta^2+t(C_3\delta^{\frac{7}{2}-2\alpha}+1)\Big\}.\label{EC3}
\end{align}
{\bf (3)} For $t\in(0,1)$ and $u\in\mR^d$, define a random set
$$
\Omega^u_t(x;\hbar):=\left\{\sup_{s\in[0,t]}|uK^\delta_s(x;\hbar)|^2\leq\delta^{-\gamma}\right\},
$$
and let
\begin{align}
\cJ^u_t(x;\hbar):=&\exp\Bigg\{c_1\delta^{\gamma-\frac{1}{2}}\int^t_0|uK^\delta_s(x;\hbar)[b,V](X^\delta_s(x;\hbar))|^2\dif s\no\\
&-\delta^{\gamma-8}\int^t_0|uK^\delta_s(x;\hbar)V(X^\delta_s(x;\hbar))|^2\dif s\Bigg\}1_{\Omega^u_t(x;\hbar)}.\label{Def0}
\end{align}
Since $\Omega^u_t(x;0)\subset\{\tau_u(x)\geq t\}$, by (\ref{EC3}) we have
\begin{align}
\sup_{x\in\mR^d}\sup_{u\in\mR^d}\mE\cJ^u_t(x;0)\leq\exp\Big\{\delta^2+t(C_3\delta^{\frac{7}{2}-2\alpha}+1)\Big\}.\label{EC4}
\end{align}
Let $0=t_0<t_1<\cdots<t_n\leq t_{n+1}=t$ be the jump times of $\hbar$. If we set
$$
\phi_{t_j}(x;\hbar):=X^\delta_{t_j-}(x;\hbar)+\sigma(X^\delta_{t_{j-}}(x;\hbar),\Delta \hbar_{t_j}),
$$
then for $s\in[0, t_{j+1}-t_j)$,
$$
X^\delta_{s+t_j}(x;\hbar)=X^\delta_s(\phi_{t_j}(x;\hbar); 0)\Rightarrow K^\delta_{s+t_j}(x;\hbar)=[\nabla\phi_{t_j}(x;\hbar)]^{-1} K^\delta_s(\phi_{t_j}(x;\hbar); 0)
$$
and
$$
\Omega^u_{t_{j+1}}(x;\hbar)=\Omega^u_{t_j}(x;\hbar)\cap\left\{\sup_{s\in[0,t_{j+1}-t_j]}|uK^\delta_{s+t_j}(x;\hbar)|^2\leq\delta^{-\gamma}\right\}.
$$
Thus, by the Markovian property,  we have for all $u\in\mR^d$,
\begin{align}
\mE\cJ^u_{t_{n+1}}(x;\hbar)&=\mE\left(\cJ^u_{t_n}(x;\hbar)
\cdot\Big(\mE\cJ^{u'}_{t_{n+1}-t_n}(y;0)\Big)\Big|_{u'=u[\nabla\phi_{t_j}(x;\hbar)]^{-1}, y=\phi_{t_j}(x;\hbar)}\right)\no\\
&\!\!\!\stackrel{(\ref{EC4})}{\leq}\mE\cJ^u_{t_n}(x;\hbar)\exp\Big\{\delta^2+(t_{n+1}-t_n)(C_3\delta^{\frac{7}{2}-2\alpha}+1)\Big\}\no\\
&\leq\cdots\cdots\cdots\cdots\cdots\cdots\no\\
&\leq\Pi_{j=0}^n\exp\Big\{\delta^2+(t_{j+1}-t_j)(C_3\delta^{\frac{7}{2}-2\alpha}+1)\Big\}\no\\
&=\exp\Big\{\delta^2(n+1)+t_{n+1}(C_3\delta^{\frac{7}{2}-2\alpha}+1)\Big\}.\label{LK1}
\end{align}
Let $N^\delta_t$ be the jump number of $\hat L^\delta_\cdot$ before time $t$, i.e.,
$$
N^\delta_t=\sum_{s\in(0,t]}1_{|\Delta\hat L^\delta_s|>0}=\int_{|z|>\delta}N((0,t],\dif z)=\sum_{s\in(0,t]}1_{|\Delta L_s|>\delta},
$$
which is a Poisson process with intensity $\int_{|z|>\delta}\nu(\dif z)=:\lambda_\delta$. If we let $m=[\log\delta^{-1}/\log 2]$, then by (\ref{EUT0}), we have
\begin{align}
\lambda_\delta&\leq\int_{|z|\geq 1}\nu(\dif z)+\sum_{k=0}^m\int_{2^k\delta\leq|z|\leq 2^{k+1}\delta}\nu(\dif z)\no\\
&\leq C+\sum_{k=0}^m(2^k\delta)^{-2}\int_{2^k\delta\leq|z|\leq 2^{k+1}\delta}|z|^2\nu(\dif z)\no\\
&\leq C+C\sum_{k=0}^m(2^k\delta)^{-2}(2^{k+1}\delta)^{2-\alpha}\no\\
&=C+C2^{2-\alpha}\sum_{k=0}^m (2^k\delta)^{-\alpha}\leq C\delta^{-\alpha}.\label{EUT00}
\end{align}
Recalling (\ref{EH2}), (\ref{Def0}) and the independence of $L^\delta$ and $\hat L^\delta$, we have for any $u\in\mR^d$,
\begin{align}
&\mE\left( \exp\left\{c_1\delta^{\gamma-\frac{1}{2}}\int^t_0|uK_s[b,V](X_s)|^2\dif s-\delta^{\gamma-8}\int^t_0|uK_sV(X_s)|^2\dif s\right\}1_{\Omega^u_t(x;\hat L^\delta)}\right)\no\\
&\qquad=\mE\Big(\mE\cJ^u_t(x;\hbar)|_{\hbar=\hat L^\delta}\Big)
=\sum_{n=0}^\infty\mE\left(\Big(\mE\cJ^u_t(x;\hbar)\Big)_{\hbar=\hat L^\delta}; N^\delta_t=n\right)\no\\
&\quad \ \ \stackrel{(\ref{LK1})}{\leq} \sum_{n=0}^\infty\exp\Big\{\delta^2(n+1)+t(C_3\delta^{\frac{7}{2}-2\alpha}+1)\Big\}\mP(N^\delta_t=n)\no\\
&\qquad=\exp\Big\{\delta^2+t(C_3\delta^{\frac{7}{2}-2\alpha}+1)\Big\}\sum_{n=0}^\infty\e^{\delta^2n}\frac{(t\lambda_\delta)^n}{n!}\e^{-t\lambda_\delta}\no\\
&\qquad=\exp\Big\{\delta^2+t(C_3\delta^{\frac{7}{2}-2\alpha}+1)+(\e^{\delta^2}-1)t\lambda_\delta\Big\}\no\\
&\quad \ \ \stackrel{(\ref{EUT00})}{\leq} \exp\Big\{2+C_4t\delta^{\frac{7}{2}-2\alpha}+C_5t\delta^{2-\alpha}\Big\},\ \ \forall t,\delta\in(0,1),\label{EH3}
\end{align}
where in the last step we have used that $\e^x-1\leq 3x$ for $x\in(0,1)$.
\\
\\
{\bf (4)} By (\ref{EH3}) and Chebyshev's inequality, we have for any $\beta\in(0,1)$,
\begin{align*}
&\mP\left\{c_1\delta^{\gamma-\frac{1}{2}}\int^t_0|uK_s[b,V](X_s)|^2\dif s-\delta^{\gamma-8}
\int^t_0|uK_sV(X_s)|^2\dif s\geq t\delta^{-\frac{\beta}{2}},\Omega^u_t(x;\hat L^\delta)\right\}\no\\
&\qquad\qquad\leq\exp\Big\{2+C_4t\delta^{\frac{7}{2}-2\alpha}+C_5t\delta^{2-\alpha}-t\delta^{-\frac{\beta}{2}}\Big\}
\end{align*}
and by (\ref{EUT1}), 
\begin{align*}
\mP([\Omega^u_t(x;\hat L^\delta)]^c)&=\mP\left\{\sup_{s\in[0,t]}|uK_s(x)|^2>\delta^{-\gamma}\right\}\\
&\leq \delta^{p\gamma}\mE\left(\sup_{s\in[0,t]}|uK_s(x)|^{2p}\right)\leq C_p|u|^{2p}\delta^{p\gamma},\ \forall p\geq 1.
\end{align*}
In particular, if $\beta\in(0\vee(4\alpha-7),1)$, then there exists a constant $\delta_0$ such that for all $\delta\in(0,\delta_0)$, $t\in(0,1)$ and $p\geq 1$,
\begin{align*}
&\sup_{|u|=1}\mP\left\{\int^t_0|uK_s[b,V](X_s)|^2\dif s\geq 
\frac{2t\delta^{\frac{1-\beta}{2}-\gamma}}{c_1},\int^t_0|uK_sV(X_s)|^2\dif s\leq t\delta^{8-\frac{\beta}{2}-\gamma}\right\}\\
&\qquad\qquad\leq\exp\{3-t\delta^{-\frac{\beta}{2}}\}+C_p\delta^p,
\end{align*}
which then gives the desired estimate by adjusting the constants and rescaling $\delta$.
\end{proof}

The reduced Malliavin matrix is defined by
\begin{align}
\hat\Sigma_t(x):=\int^t_0K_s(x)[(\nabla_z\sigma)(\nabla_z\sigma)^*](X_s(x),0)K^*_s(x)\dif s.\label{Ma}
\end{align}
We are now in a position to prove the following main result of this section.
\bt\label{Th0}
Let $B_1(x):=\nabla_z\sigma(x,0)$ and define $B_{j+1}(x):=b(x)\cdot\nabla B_j(x)-B_j(x)\cdot\nabla b(x)$ for $j\in\mN$.
Assume that for some $j_0\in\mN$,
\begin{align}
\inf_{x\in\mR^d}\inf_{|u|=1}\sum_{j=1}^{j_0}|uB_j(x)|^2=:c_0>0.\label{UH0}
\end{align}
Under (\ref{33}), (\ref{34}) and (\ref{Con3}), there exist $\gamma=\gamma(\alpha, j_0)\in(0,1)$ and constants $C_2\geq 1$, $c_2\in(0,1)$ 
such that for all $t\in(0,1)$, $\lambda\geq 1$ and $p\geq 1$,
\begin{align}
\sup_{|u|=1}\sup_{x\in\mR^d}\mE\exp\left\{-\lambda u\hat\Sigma_t(x)u^*\right\}\leq C_2\exp\{- c_2 t\lambda^\gamma\}+C_p(\lambda t)^{-p}.\label{EE8}
\end{align}
\et
\begin{proof}  Let $\beta_1,\beta_2,\beta_3$ be as in  (\ref{EE4}). Set $a:=\frac{\beta_1}{\beta_2}\leq 1$ and define for $j=1,\cdots, j_0$, 
$$
E_j:=\left\{\int^t_0|uK_sB_j(X_s)|^2\dif s\leq t\delta^{a^j\beta_2}\right\}.
$$
Since $a^{j+1}\beta_2=a^j\beta_1$ and $B_{j+1}=[b,B_j]$, by (\ref{EE4}) with $\delta$ replaced by $\delta^{a^j}$, we have for any $p\geq 1$,
\begin{align}
\mP(E_j\cap E_{j+1}^c)\leq C_1\exp\{-c_1t\delta^{-a^j\beta_3}\}+C_p\delta^{a^jp}.\label{EC5}
\end{align}
Noticing that
$$
E_1\subset\left(\cap_{j=1}^{j_0} E_j\right)\cup\left(\cup_{j=1}^{j_0-1}(E_j\cap E_{j+1}^c)\right),
$$
we have
\begin{align}
\mP(E_1)\leq \mP\left(\cap_{j=1}^{j_0} E_j\right)+\sum_{j=1}^{j_0-1}\mP(E_j\cap E_{j+1}^c).\label{EC6}
\end{align}
On the other hand, if we define
$$
\tau:=\inf\{t\geq 0: |J_t|\geq\delta^{-a^{j_0}\beta_2/3}\},
$$
then for any $s\leq\tau$ and $|u|=1$,
$$
|uK_s|^2\geq|J_s|^{-2}\geq\delta^{2a^{j_0}\beta_2/3}.
$$
Thus, by (\ref{UH0}) we have
\begin{align}
\bigcap_{j=1}^{j_0}E_j\cap\{\tau\geq t\}&\subset\left\{\sum_{j=1}^{j_0}\int^t_0|uK_sB_j(X_s)|^2\dif s\leq t\sum_{j=1}^{j_0}\delta^{a^j\beta_2},\tau\geq t\right\}\no\\
&\subset\left\{c_0\int^t_0|uK_s|^2\dif s\leq t \sum_{j=1}^{j_0}\delta^{a^{j_0}\beta_2},\tau\geq t\right\}\no\\
&\subset\left\{tc_0\delta^{2a^{j_0}\beta_2/3}\leq t j_0\delta^{a^{j_0}\beta_2}\right\}=\emptyset,\label{EC7}
\end{align}
provided $\delta<\delta_1=(c_0/j_0)^{3/(a^{j_0}\beta_2)}$. On the other hand, by (\ref{EUT1}), we have for any $p\geq 2$,
$$
\mP(\tau<t)\leq \mP\left(\sup_{s\in[0,t]}|J_s|\geq\delta^{-a^{j_0}\beta_2/3}\right)\leq \mE\left(\sup_{s\in[0,t]}|J_s|^p\right)\delta^{pa^{j_0}\beta_2/3}
\leq C_p\delta^{pa^{j_0}\beta_2/3}.
$$
Therefore, combining (\ref{EC5})-(\ref{EC7}) and resetting $\eps=\delta^{\beta_1}$ and $\theta=a^{j_0}\beta_3/\beta_1$, we obtain that for all $\eps\in(0,1)$, $t\in(0,1)$ and $p\geq 1$,
$$
\sup_{|u|=1}\mP\left\{\int^t_0|uK_sB_1(X_s)|^2\dif s\leq t\eps\right\}\leq C_2\exp\Big\{-c_1t\eps^{-\theta}\Big\}+C_p\eps^p.
$$
For $\lambda\geq t$, setting $r:=(\lambda/t)^{\frac{-1}{1+\theta}}$ and $\xi:=\frac{1}{t}\int^t_0|uK_sB_1(X_s)|^2\dif s$, we have
\begin{align*}
\mE\e^{-\lambda\xi}&=\int^\infty_0\lambda\e^{-\lambda\eps}\mP(\xi\leq\eps)\dif \eps\\
&\leq\int^\infty_r\lambda\e^{-\lambda\eps}\dif\eps+C\int^r_0\lambda\e^{-\lambda\eps}(\e^{-c_1t\eps^{-\theta}}+\eps^p)\dif \eps\\
&=\e^{-\lambda r}+C\int^{\lambda r}_0\e^{-s-c_1t \lambda^\theta s^{-\theta}}\dif s+C\lambda^{-p}\int^{\lambda r}_0\e^{-s}s^p\dif s\\
&\leq\e^{-\lambda r}+C\e^{-c_1t r^{-\theta}}\int^{\lambda r}_0\e^{-s}\dif s+C\lambda^{-p}\\
&\leq\e^{-t (\lambda/t)^{\frac{\theta}{1+\theta}}}+C\e^{-c_1t (\lambda/t)^{\frac{\theta}{1+\theta}}}+C\lambda^{-p}.
\end{align*}
By replacing $\lambda$ with $\lambda t$ and recalling (\ref{Ma}), we obtain the desired estimate (\ref{EE8}).
\end{proof}

\section{Smooth densities for SDEs without big jumps}
In the remainder of this paper, we assume {\bf (H$^\sigma_b$)} and {\bf (H$^\nu$)} and choose $\delta$ being small enough so that 
\begin{align}
|\nabla_x\sigma(x,z)|\leq \tfrac{1}{2},\ \ |z|\leq\delta,\label{Con}
\end{align}
and set
$$
\Gamma^\delta_0:=\{z\in\mR^d: 0<|z|<\delta\}.
$$
Let $X_t(x)=X_t$ solve the following SDE:
\begin{align}
X_t=x+\int^t_0b(X_s)\dif s+\int^t_0\!\!\!\int_{\Gamma^\delta_0}\sigma(X_{s-},z)\tilde N(\dif s,\dif z).\label{SDE9}
\end{align}
It is well known that the generator of $X_t(x)$ is given by $\cL_0$ in (\ref{Op3}).

This section is based on Subsection 2.1, Lemma \ref{Le25} and Theorem \ref{Th0}. We first prove the following Malliavin differentiability of $X_t$ with respect to $\omega$ in the sense of Theorem \ref{Th21}.
\bl
Fix $\v\in\mV_{\infty-}$. For any $t\in[0,1]$, we have $X_t\in\mW^{1,\infty-}_\v(\Omega)$ and
\begin{align}
D_\v X_t&=\int^t_0\nabla b(X_s)D_\v X_s\dif s+\int^t_0\!\!\!\int_{\Gamma^\delta_0}\nabla_x \sigma(X_{s-},z)D_\v X_{s-}\tilde N(\dif s,\dif z)\no\\
&\quad+\int^t_0\!\!\!\int_{\Gamma^\delta_0}\<\nabla_z\sigma(X_{s-},z),\v(s,z)\> N(\dif s,\dif z).\label{SDE7}
\end{align}
Moreover, for any $p\geq 2$, we have
\begin{align}
\sup_{x\in\mR^d}\mE\left(\sup_{t\in[0,1]}|D_\v X_t(x)|^p\right)<\infty.\label{HH2}
\end{align}
\el
\begin{proof}
{\bf (1)} Consider the following Picard's iteration: $X_t^0\equiv x$ and for $n\in\mN$,
\begin{align*}
X^{n}_{t}:=x+\int^t_0b(X^{n-1}_{s})\dif s+\int^t_0\!\!\!\int_{\Gamma^\delta_0}\sigma_2(X^{n-1}_{s-},z)\tilde N(\dif s,\dif z).
\end{align*}
Since $b$ and $\sigma$ are Lipschitz continuous, it is by now standard to prove that for any $p\geq 2$,
\begin{align}
\sup_{n\in\mN}\mE\left(\sup_{t\in[0,1]}|X_t^n|^p\right)<\infty\mbox{ and }\lim_{n\rightarrow \infty}\mE\left(\sup_{t\in[0,1]}|X_t^n-X_t|^p\right)=0.\label{pp}
\end{align}
{\bf (2)} Now we use the induction to prove that for each $n\in\mN$,
\begin{align}
X^n_t\in\mW^{1,\infty-}_\v(\Omega)\ \mbox{ and }\ \ \mE\left(\sup_{t\in[0,1]}|D_\v X^{n}_t|^p\right)<+\infty,\ \ \forall p\geq 2.\label{EL1}
\end{align}
First of all, it is clear that (\ref{EL1}) holds for $n=0$. Suppose now that (\ref{EL1}) holds for some $n\in\mN$.
By (\ref{pp}) and the induction hypothesis, it is easy to check that the assumptions of Proposition \ref{Pr1} are satisfied.
Thus, $X^{n+1}_t\in\mW^{1,\infty-}_\v(\Omega)$ and
\begin{align*}
D_\v X^{n+1}_t&=\int^t_0\nabla b(X^n_{s})D_\v X^n_{s}\dif s+\int^t_0\!\!\!\int_{\Gamma^\delta_0}\nabla_x\sigma(X^n_{s-},z)D_\v X^n_{s-}\tilde N(\dif s,\dif z)\no\\
&\quad+\int^t_0\!\!\!\int_{\Gamma^\delta_0}\<\nabla_z\sigma(X^n_{s-},z),\v(s,z)\> N(\dif s,\dif z).
\end{align*}
By Lemma \ref{Le2}, we have for any $p\geq 2$,
\begin{align*}
\mE\left(\sup_{s\in[0,t]}|D_\v X^{n+1}_s|^p\right)&\leq C\int^t_0\mE|D_\v X^n_{s}|^p\dif s
+C\mE\left(\int^t_0\!\!\!\int_{\Gamma^\delta_0}|\<\nabla_z\sigma(X^n_{s-},z),\v(s,z)\>|\nu(\dif z)\dif s\right)^p\no\\
&\quad+C\mE\left(\int^t_0\!\!\!\int_{\Gamma^\delta_0}|\<\nabla_z\sigma(X^n_{s-},z),\v(s,z)\>|^p\nu(\dif z)\dif s\right).
\end{align*}
Since $\v\in\mV_{\infty-}$, by {\bf (H$^\sigma_b$)} we further have
\begin{align*}
\mE\left(\sup_{s\in[0,t]}|D_\v X^{n+1}_s|^p\right)&\leq C\int^t_0\mE|D_\v X^n_{s}|^p\dif s+C
\leq C\int^t_0\mE\left(\sup_{r\in[0,s]}|D_\v X^n_r|^p\right)\dif s+C,
\end{align*}
where $C$ is independent of $n$ and the starting point $x$. Thus, we have proved (\ref{EL1}) by the induction
hypothesis. Moreover, by Gronwall's inequality, we also have
\begin{align}
\sup_{n\in\mN}\mE\left(\sup_{s\in[0,1]}|D_\v X^n_s|^p\right)<+\infty.\label{HH}
\end{align}
{\bf (3)} Let $Y_t$ solve the following linear matrix-valued SDE:
\begin{align*}
Y_t&=\int^t_0\nabla b(X_s)Y_s\dif s+\int^t_0\!\!\!\int_{\Gamma^\delta_0}\nabla_x\sigma(X_{s-},z)Y_{s-}\tilde N(\dif s,\dif z)\no\\
&\qquad+\int^t_0\!\!\!\int_{\Gamma^\delta_0}\<\nabla_z\sigma(X_{s-},z),\v(s,z)\> N(\dif s,\dif z).
\end{align*}
By Fatou's lemma and (\ref{pp}), (\ref{HH}), for any $p\geq 2$, we have
\begin{align*}
\varlimsup_{n\to\infty}\mE|D_\v X^n_t-Y_t|^p&\leq C\int^t_0\varlimsup_{n\to\infty}\mE|D_\v X^{n-1}_s-Y_s|^p\dif s,
\end{align*}
which then gives
$$
\varlimsup_{n\to\infty}\mE|D_\v X^n_t-Y_t|^p=0.
$$
Thus, $X_t\in\mW^{1,p}_\v(\Omega)$ and $D_\v X_t=Y_t$. Moreover, the estimate (\ref{HH2}) follows by (\ref{pp}) and (\ref{HH}).
\end{proof}
Let $J_t=J_t(x)$ be the Jacobian matrix of $x\mapsto X_t(x)$, and $K_t(x)$ the inverse of $J_t(x)$. 
Recalling equations (\ref{SDE2}) and (\ref{SDE7}), by the formula of constant variation, we have for any $\v\in\mV_{\infty-}$,
\begin{align}
D_\v X_t=J_t\int^t_0\!\!\!\int_{\Gamma^\delta_0}K_s\nabla_z\sigma(X_{s-},z)\v(s,z) N(\dif s,\dif z).\label{EQ1}
\end{align}
Here the integral is the Lebesgue-Stieltjes integral. Let
$$
U(x,z):=(\mI+\nabla_x\sigma(x,z))^{-1}\nabla_z\sigma(x,z),\ \ x\in\mR^d,\ \ z\in\Gamma^\delta_0,
$$
and define
$$
\v_j(x;s,z):=[K_{s-}(x)U(X_{s-}(x),z)]^*_{\cdot j}\zeta(z),
$$
where $\zeta(z)=\zeta_\delta(z)$ is a nonnegative smooth function with
$$
\zeta_\delta(z)=|z|^3,\ \ |z|\leq \delta/4,\ \ \zeta_\delta(z)=0,\ \ |z|>\delta/2.
$$
The following lemma is easy to be verified by definitions and (\ref{EUT1}).
\bl
For any $m\in\mN_0$, there is a constant $C>0$ such that for all $x\in\mR^d$ and $z\in\Gamma^\delta_0$,
\begin{align}
|\nabla^m_xU(x,z)|,\ |\nabla^m_zU(x,z)|\leq C,\ \ |U(x,z)-U(x,0)|\leq C|z|.\label{TYR1}
\end{align}
Moreover, for each $j=1,\cdots, d$ and $x\in\mR^d$, $\v_j(x)\in\mV_{\infty-}$.
\el
Write
$$
\Theta(s,z):=\Theta(x;s,z):=(\v_1(x;s,z),\cdots,\v_d(x;s,z))
$$
and
$$
(D_\Theta X_t)_{ij}:=D_{\v_j}X^i_t.
$$
Since by equation (\ref{SDE3}),
$$
K_s=K_{s-}(\mI+\nabla_x\sigma(X_{s-},\Delta L_s))^{-1},
$$
by (\ref{EQ1}) we have
\begin{align}
D_\Theta X_t(x)=J_t(x)\Sigma_t(x),\label{Var}
\end{align}
where
\begin{align}
\Sigma_t(x):=\int^t_0\!\!\!\int_{\Gamma^\delta_0}K_{s-}(x)(UU^*)(X_{s-}(x),z)K^*_{s-}(x)\zeta(z)N(\dif s,\dif z).\label{Var1}
\end{align}

\bl\label{Le43}
For any $p\geq 2$ and $m,k\in\mN_0$ with $m+k\geq 1$, we have
\begin{align}
&\sup_{x\in\mR^d}\mE\left(\sup_{t\in[0,1]}|D_{\v_{j_1}}\cdots D_{\v_{j_m}}\nabla^k X_t(x)|^p\right)<\infty,\label{TR4}\\
&\sup_{x\in\mR^d}\mE\left(\sup_{t\in[0,1]}|D_{\v_{j_1}}\cdots D_{\v_{j_m}}\div(\v_i(x))|^p\right)<\infty,\label{TR5}
\end{align}
where $j_1,\cdots,j_m$ and $i$ runs in $\{1,2,\cdots,d\}$.
\el
\begin{proof}
For $m+k=1$, (\ref{TR4}) has been proven in (\ref{EUT1}) and (\ref{HH2}). For general $k$ and $m$, it follows by induction.
Let us look at (\ref{TR5}) with $m=1$. Notice that by (\ref{ER1}),
\begin{align*}
\div(\v_i)=\int^1_0\!\!\!\int_{\Gamma^\delta_0}[\<\nabla\log\kappa(z),\v_i(s,z)\>+\div_z(\v_i)(s,z)]\tilde N(\dif s,\dif z).
\end{align*}
By Proposition \ref{Pr1}, we have
\begin{align*}
&D_{\v_j}\div(\v_i)=\int^1_0\!\!\!\int_{\Gamma^\delta_0}[\<\nabla\log\kappa(z),D_{\v_j}\v_i(s,z)\>+D_{\v_j}\div_z(\v_i)(s,z)]\tilde N(\dif s,\dif z)\\
&\quad+\int^1_0\!\!\!\int_{\Gamma^\delta_0}\<\v_j(s,z),\nabla_z\<\nabla\log\kappa(z),\v_i(s,z)\>+\nabla_z\div_z(\v_i)(s,z)\>N(\dif s,\dif z).
\end{align*}
In view of $\mathrm{supp}{\v_i}(s,\cdot)\subset\Gamma^{\frac{\delta}{2}}_0$,
by Lemma \ref{Le2} and (\ref{ET01}), (\ref{TYR1}), (\ref{TR4}), one obtains (\ref{TR5}) with $m=1$. For general $m$, it follows by similar calculations.
\end{proof}

Below we define
\begin{align}
\cT^0_t f(x):=\mE f(X_t(x)).\label{Semi}
\end{align}
The following lemma is proven in appendix.
\bl\label{Le41}
Under {\bf (H$^\sigma_b$)}, there exists a constant $C>0$ such that for any $f\in L^1(\mR^d)$,
\begin{align}
\sup_{t\in[0,1]}\int_{\mR^d}|\cT^0_tf(x)|\dif x\leq C\int_{\mR^d}|f(x)|\dif x.\label{IU7}
\end{align}
\el 

Now we can prove the following main result of this section.
\bt\label{Th45}
Assume {\bf (H$^\sigma_b$), (H$^\nu$)} and {\bf (UH)} and let $\delta$ be as in (\ref{Con}).
For any $k,m,n\in\mN_0$ and $p\in(1,\infty]$, there exist $\gamma_{kmn}\geq 0$ only depending on $k,m,n,\alpha,j_0,d$ and
a constant $C\geq 1$ such that for all $f\in \mW^{n,p}(\mR^d)$ and $t\in(0,1)$,
\begin{align}
\|\nabla^k\cT^0_t\nabla^{m}f\|_p\leq C t^{-\gamma_{kmn}}\|f\|_{n,p},\label{NJ6}
\end{align}
where $\gamma_{kmn}$ is increasing with respect to $k,m$ and decreasing in $n$, and $\gamma_{kmn}=0$ for $n\geq k+m$.
In particular, $X_t(x)$ admits a smooth density $\rho_t(x,y)$ with $\rho_t\in C^\infty_b(\mR^d\times\mR^d)$ such that
$$
\p_t\rho_t(x,y)=\cL_0\rho_t(\cdot,y)(x),\ \ \forall(t,x,y)\in(0,1)\times\mR^d\times\mR^d.
$$
\et
\begin{proof}
Below we only prove (\ref{NJ6}) for $p\in(1,\infty)$. For $p=\infty$, it is similar and simpler.
We assume $f\in C^\infty_0(\mR^d)$ and divide the proof into four steps.
\\
\\
{\bf (1)} Let $\Sigma_t(x)$ be defined by (\ref{Var1}). In view of $U(x,0)=\nabla_z\sigma(x,0)$, by (\ref{TR7}) and (\ref{EE8}), there are constants $C_3\geq 1$, $c_3\in(0,1)$ and
$\gamma=\gamma(\alpha,j_0)\in(0,1)$ such that for all $t\in(0,1)$, $\lambda\geq 1$ and $p\geq 1$,
\begin{align}
\sup_{|u|=1}\sup_{x\in\mR^d}\mE\exp\left\{-\lambda u\Sigma_t(x)u^*\right\}\leq C_3\exp\{- c_3 t\lambda^\gamma\}+C_p(\lambda t)^{-p}.\label{EE88}
\end{align}
where $\Sigma_t(x)$ is defined by (\ref{Var1}).
As in \cite[Lemma 5.3]{Zh4}, for any $p\geq 1$, there exist constant $C\geq 1$ an $\gamma'=\gamma'(\alpha,j_0,d)$ such that for all $t\in(0,1)$,
$$
\sup_{x\in\mR^d}\mE \Big((\det\Sigma_t(x))^{-p}\Big)\leq Ct^{-\gamma' p},
$$
which in turn gives that for all $p\geq 1$,
\begin{align}
\sup_{x\in\mR^d}\|\Sigma^{-1}_t(x)\|_{L^p(\Omega)}\leq Ct^{-\gamma'}.\label{EK9}
\end{align}
{\bf (2)} For $t\in(0,1)$ and $x\in\mR^d$, let $\sC_t(x)$ be the class of all polynomial functionals of 
$$
\div\Theta, \Sigma^{-1}_t, K_t, \big(\nabla^k X_t\big)_{k=1}^{\ell_1}, \big(D_{\v_{j_1}}\cdots D_{\v_{j_m}}(X_t,\cdots,\nabla^{\ell_2} X_t, K_t, \div\Theta, \Sigma_t)\big)_{m=1}^{\ell_3},
$$
where $\ell_1,\ell_2,\ell_3\in\mN$, $j_i\in\{1,\cdots,d\}$ and the starting point $x$ is dropped in the above random variables. 
By (\ref{EK9}) and Lemma \ref{Le43}, for any $H_t(x)\in \sC_t(x)$, there exists a $\gamma(H)\geq 0$ only depending on the degree of $\Sigma^{-1}_t$ and
$\alpha,j_0,d$ such that for all $t\in(0,1)$ and $p\geq 1$,
\begin{align}
\sup_{x\in\mR^d}\|H_t(x)\|_{L^p(\Omega)}\leq C_p t^{-\gamma(H)}.\label{HHJ}
\end{align}
Notice that if $H_t$ does not contain $\Sigma_t^{-1}$, then $\gamma(H)=0$.
\\
\\
{\bf (3)} Let $\xi\in\sC_t(x)$. Recalling that $D_\Theta X$ is an invertible matrix, by (\ref{Var}) and the integration by parts formula (\ref{ER88}), we have
\begin{align*}
\mE\Big((\nabla f)(X_t)\xi\Big)
&=\mE\Big(\nabla f(X_t)D_\Theta X_t\cdot (D_\Theta X_t)^{-1}\xi\Big)=\mE\Big(D_\Theta f(X_t)\Sigma_t^{-1}K_t\xi\Big)\\
&=\mE\Big(-f(X_t)(\div\Theta\cdot\Sigma_t^{-1}K_t\xi-D_{\v_i}[(\Sigma_t^{-1}K_t)_{i\cdot}\xi])\Big)=\mE(f(X_t)\xi'),
\end{align*}
where $\xi'\in\sC_t(x)$. Starting from this formula, by the chain rule and induction, we have
\begin{align*}
&\nabla^k\mE\Big((\nabla^mf)(X_t)\Big)
=\sum_{j=0}^k\mE\Big((\nabla^{m+j}f)(X_t)G_j(\nabla X_t,\cdots,\nabla^k X_t)\Big)=\sum_{j=0}^n\mE \Big((\nabla^jf)(X_t)H_j\Big),
\end{align*}
where $\{G_j, j=1,\cdots,k\}$ are real polynomial functions and $H_j\in\sC_t(x)$. Notice that if $n=k+m$, then $H_j$ will not contain $\Sigma^{-1}_t$.
\\
\\
{\bf (4)} Now,  for any $p\in(1,\infty)$, by H\"older's inequality, we have
\begin{align*}
\|\nabla^k\cT^0_t\nabla^{m}f\|_{L^p(\mR^d)}&\leq \sum_{j=0}^n\left(\int_{\mR^d}\left|\mE \Big((\nabla^jf)(X_t(x))H_j(x)\Big)\right|^p\dif x\right)^{\frac{1}{p}}\\
&\leq\sum_{j=0}^n\left(\int_{\mR^d}\mE \Big(|\nabla^jf|^p(X_t(x))\Big)\Big(\mE|H_j(x)|^{\frac{p}{p-1}}\Big)^{p-1}\dif x\right)^{\frac{1}{p}}\\
&\stackrel{(\ref{HHJ})}{\leq} C\sum_{j=0}^n t^{-\gamma(H_j)}\left(\int_{\mR^d}\mE \Big(|\nabla^jf|^p(X_t(x))\Big)\dif x\right)^{\frac{1}{p}}\\
&\stackrel{(\ref{IU7})}{\leq} Ct^{-\max\{\gamma(H_j), j=1,\cdots, n\}}\|f\|_{n,p},\ \ t\in(0,1).
\end{align*}
As for the second conclusion, it follows by (\ref{NJ6}) and Sobolev's embedding theorem (cf. \cite{Nu}). The proof is thus complete.
\end{proof}

\section{Proofs of Theorems \ref{Main} and \ref{Main1}}
We first recall some definitions about the Sobolev and H\"older spaces.
For $k\in\mN_0$ and $p\in[1,\infty]$, let $\mW^{k,p}=\mW^{k,p}(\mR^d)$ be the usual Sobolev spaces with the norm:
$$
\|\varphi\|_{k,p}:=\sum_{j=0}^k\|\nabla^j \varphi\|_p.
$$
For $\beta\geq 0$ and $p\in[1,\infty)$, let $\mH^{\beta,p}:=(I-\Delta)^{-\frac{\beta}{2}}(L^p(\mR^d))$ be the usual Bessel potential space. 
For $p=\infty$, let $\mH^{\beta,\infty}$ be the usual H\"older space, i.e., if $\beta=k+\theta$ with $\theta\in[0,1)$, then
$$
\|\varphi\|_{\beta,\infty}:=\|\varphi\|_{k,\infty}+[\nabla^k\varphi]_{\theta}<\infty,
$$
where $[\nabla^k\varphi]_0:=0$ by convention and for $\theta\in(0,1)$,
$$
[\nabla^k\varphi]_\theta:=\sup_{x\not=y}\frac{|\nabla^k\varphi(x)-\nabla^k\varphi(y)|}{|x-y|^{\theta}}.
$$
Notice that $\mH^{k,\infty}=C^k_b(\mR^d)$ for $k\in\mN_0$.
It is well known that for any $k\in\mN_0$ and $p\in(1,\infty)$ (cf. \cite{St}), 
$$
\mH^{k,p}=\mW^{k,p},
$$ 
and for any $\beta_1,\beta_2\geq 0$, $p\in(1,\infty)$ and $\theta\in[0,1]$,
\begin{align}
[\mH^{\beta_1,p}, \mH^{\beta_2,p}]_\theta=\mH^{\beta_1+\theta(\beta_2-\beta_1),p},\label{Int1}
\end{align}
and if $\beta_1+\theta(\beta_2-\beta_1)$ is not an integer, then
\begin{align}
(\mH^{\beta_1,\infty}, \mH^{\beta_2,\infty})_{\theta,\infty}=\mH^{\beta_1+\theta(\beta_2-\beta_1),\infty},\label{Int2}
\end{align}
where $[\cdot,\cdot]_\theta$ (resp. $(\cdot,\cdot)_{\theta,\infty}$) stands for the complex (resp. real) interpolation space.

We recall the following interpolation theorem (cf. \cite[p.59, Theorem (a)]{Tr}).
\bt\label{Th51}
Let $A_i\subset B_i, i=0,1$ be Banach spaces. Let $\sT:A_i\to B_i, i=0,1$ be bounded linear operators.
For $\theta\in[0,1]$,  we have
$$
\|\sT\|_{A_\theta\to B_\theta}\leq\|\sT\|^{1-\theta}_{A_0\to B_0}\|\sT\|^{\theta}_{A_1\to B_1},
$$
where $A_\theta:=[A_0,A_1]_{\theta}$, $B_\theta:=[B_0,B_1]_{\theta}$, and
$\|\sT\|_{A_\theta\to B_\theta}$ denotes the operator norm of $\sT$ mapping $A_\theta$ to $B_\theta$.
The same is true for real interpolation spaces.
\et

Let $\cT^0_t$ be the semigroup defined by (\ref{Semi}), whose generator is given by $\cL_0$. We have
\bl\label{Le53}
Let $\gamma_{100}$ be the same as in Theorem \ref{Th45}.
For any $p\in(1,\infty)$, $\theta\in[0,1)$ and $\beta\geq 0$, there exit constants $C_1, C_2>0$ such that for all $t\in(0,1)$,
\begin{align}
\|\cT^0_t\varphi\|_{\theta+\beta,p}\leq C_1t^{-\theta\gamma_{100}}\|\varphi\|_{\beta,p},\label{HG2}
\end{align}
and if $\beta$ and $\theta+\beta$ are not integers, then
\begin{align}
\|\cT^0_t\varphi\|_{\theta+\beta,\infty}\leq C_2t^{-\theta\gamma_{100}}\|\varphi\|_{\beta,\infty}.\label{HG3}
\end{align}
\el
\begin{proof}
Let $\theta\in[0,1)$ and $\beta\geq 0$. For any $p\in(1,\infty]$,
by Theorem \ref{Th45} and interpolation Theorem \ref{Th51}, there exists a constant $C>0$ such that for all $t\in(0,1)$,
$$
\|\cT^0_t\varphi\|_{\frac{\beta}{1-\theta},p}\leq C\|\varphi\|_{\frac{\beta}{1-\theta},p}
$$
and
$$
\|\cT^0_t\varphi\|_{1,p}\leq Ct^{-\gamma_{100}}\|\varphi\|_{p}.
$$
On the other hand, noticing that by (\ref{Int1}),
$$
[\mH^{\frac{\beta}{1-\theta},p},\mH^{1,p}]_\theta=\mH^{\beta+\theta,p},\ \
[\mH^{\frac{\beta}{1-\theta},p},\mH^{0,p}]_\theta=\mH^{\beta,p},
$$
and if $\beta$ and $\theta+\beta$ are not integers, then by (\ref{Int2}),
$$
(\mH^{\frac{\beta}{1-\theta},\infty},\mH^{1,\infty})_{\theta,\infty}=\mH^{\beta+\theta,\infty},\ \
(\mH^{\frac{\beta}{1-\theta},\infty},\mH^{0,\infty})_{\theta,\infty}=\mH^{\beta,\infty},
$$
by interpolation Theorem \ref{Th51} again, we obtain the desired estimate.
\end{proof}
\subsection{Proof of Theorem \ref{Main}}
Let $\sL$ be a bounded linear operator in $C_b(\mR^d)$ and Sobolev spaces
$\mW^{k,p}(\mR^d)$ for any $p>1$ and $k\in\mN_0$. Let $\cT_t$ be the semigroup in $L^p(\mR^d)$ associated with $\cL_0+\sL$, i.e.,
for any $\varphi\in L^p(\mR^d)$,
$$
\p_t\cT_t\varphi=\cL_0\cT_t\varphi+\sL\cT_t\varphi.
$$
By Duhamel's formula, we have
\begin{align}
\cT_t\varphi=\cT^0_t\varphi+\int^t_0\cT^0_{t-s}\sL\cT_s\varphi\dif s.\label{NJ55}
\end{align}
\bl\label{Le52}
Let $\gamma_{100}$ be as in Theorem \ref{Th45}. Fix $\theta\in(0,\frac{1}{\gamma_{100}}\wedge 1)$. For any $m\in\mN$ and $p\in(1,\infty)$, there exists a constant $C>0$ such that for all $t\in(0,1)$ and $\varphi\in L^p(\mR^d)$,
\begin{align}
\|\cT_t\varphi\|_{m\theta,p}\leq Ct^{-m\theta\gamma_{100}}\|\varphi\|_p.\label{KJ1}
\end{align}
\el
\begin{proof}
First of all, since $\sL$ is a bounded linear operator in $\mW^{k,p}$, by interpolation Theorem \ref{Th51}, we have for all $\beta\geq 0$ and $p\in(1,\infty)$,
$$
\|\sL\varphi\|_{\beta,p}\leq C\|\varphi\|_{\beta, p}.
$$
Let $\theta\in(0,\frac{1}{\gamma_{100}}\wedge 1)$ and $m\in\mN$. By (\ref{NJ55}) and Lemma \ref{Le53}, we have
\begin{align*}
\|\cT_t\varphi\|_{m\theta,p}&\leq\|\cT^0_t\varphi\|_{m\theta,p}+\int^t_0\|\cT^0_{t-s}\sL\cT_s\varphi\|_{m\theta,p}\dif s\\
&\leq Ct^{-\theta\gamma_{100}}\|\varphi\|_{(m-1)\theta,p}+C\int^t_0\|\cT_s\varphi\|_{m\theta,p}\dif s,
\end{align*}
which, by Gronwall's inequality, yields that for  all $t\in(0,1)$,
$$
\|\cT_t\varphi\|_{m\theta,p}\leq Ct^{-\theta\gamma_{100}}\|\varphi\|_{(m-1)\theta, p}.
$$
Thus, by iteration we obtain
\begin{align*}
\|\cT_{mt}\varphi\|_{m\theta,p}\leq Ct^{-\theta\gamma_{100}}\|\cT_{(m-1)t}\varphi\|_{(m-1)\theta,p}\leq\cdots
\leq Ct^{-m\theta\gamma_{100}}\|\varphi\|_{p},
\end{align*}
which gives the estimate (\ref{KJ1}) by resetting $mt$ with $t$.
\end{proof}

Now we can give
\begin{proof}[Proof of Theorem \ref{Main}]
For any $p\in(1,\infty)$ and $\varphi\in L^p(\mR^d)$, by Lemma \ref{Le52}
and Sobolev's embedding theorem, we have $\cT_t\varphi\in C^\infty_b(\mR^d)$ and for any $k\in\mN_0$ and $t\in(0,1)$,
\begin{align}
\|\cT_t\varphi\|_{k,\infty}\leq C\|\cT_t\varphi\|_{k+d,p}\leq Ct^{-(k+d)\gamma_{100}}\|\varphi\|_p.\label{KG1}
\end{align}
In particular, there is a function $\rho_t(x,\cdot)\in L^{\frac{p}{p-1}}(\mR^d)$ such that for any $\varphi\in L^p(\mR^d)$,
$$
\cT_t\varphi(x)=\int_{\mR^d}\varphi(y)\rho_t(x,y)\dif y.
$$
By (\ref{KG1}), we obtain
$$
\|\nabla^k_x\rho_t(x,\cdot)\|_{\frac{p}{p-1}}\leq Ct^{-(k+d)\gamma_{100}},
$$
where $\nabla^k_x$ stands for the distributional derivative, $C$ is independent of $x$. By Fubini's theorem, we have for any $R>0$,
$$
\int_{\mR^d}\!\int_{B_R}|\nabla^k_x\rho_t(x,y)|^\frac{p}{p-1}\dif x\dif y<\infty,
$$
which, by Sobolev's embedding theorem again, produces that for almost all $y\in\mR^d$,
$$
x\mapsto\rho_t(x,y)\mbox{ is smooth}.
$$
As for equation (\ref{NG1}), it follows by (\ref{NJ55}).
\end{proof}
\subsection{Proof of Theorem \ref{Main1}}
Let $\delta$ be as in (\ref{Con}). We decompose the operator $\cL^\nu_{\sigma,b}$ as 
$$
\cL^\nu_{\sigma,b}f(x)=\cL_0f(x)+\cL_1f(x),
$$
where
$$
\cL_0 f(x):=\mathrm{p.v.}\int_{|z|<\delta}f(x+\sigma(x,z))-f(x))\nu(\dif z)+b(x)\cdot\nabla f(x)
$$
and
$$
\cL_1 f(x):=\int_{|z|\geq \delta}f(x+\sigma(x,z))-f(x))\nu(\dif z).
$$
\bl
If $\int_{|z|\geq 1}|z|^q\nu(\dif z)<\infty$ for some $q>0$, then for any $\beta\in[0,q]$,
there exists a constant $C>0$ such that for all $f\in\mH^{\beta,\infty}$,
\begin{align}
\|\cL_1 f\|_{\beta,\infty}\leq C\|f\|_{\beta,\infty}.\label{ENL}
\end{align}
\el
\begin{proof}
First of all, (\ref{ENL}) is clearly true for $\beta=0$.
By interpolation Theorem \ref{Th51}, it suffices to prove (\ref{ENL}) for $\beta\in[0,q]\cap\mN$ and $\beta=q$. 
Setting $\phi_z(x):=x+\sigma(x,z)$, by (\ref{Con1}), we have
\begin{align}
\|\nabla^m_x\phi_z\|_\infty\leq C(1+|z|),\ \forall m\in\mN.\label{JH1}
\end{align}
If $q\in(0,1)$,  then
\begin{align*}
[f\circ\phi_z]_q&=\sup_{x\not= y}\frac{|f\circ\phi_z(x)- f\circ\phi_z(y)|}{|x-y|^q}
\leq[f]_q\sup_{x\not= y}\frac{|\phi_z(x)- \phi_z(y)|^q}{|x-y|^q}\\
&\leq[f]_q\|\nabla\phi_z\|_\infty^q\stackrel{(\ref{JH1})}{\leq} C[f]_q(1+|z|^q).
\end{align*}
Hence,
$$
[\cL_1 f]_q\leq C[f]_q\int_{|z|\geq 1}(1+|z|^q)\nu(\dif z).
$$
For $q=1$, it is easy to see that (\ref{ENL}) is true by the chain rule. Now assume $q\in(1,2)$.
By the chain rule we have
\begin{align*}
[\nabla (f\circ\phi_z)]_{q-1}
&=\sup_{x\not= y}\frac{|(\nabla f)\circ\phi_z(x)\cdot\nabla\phi_z(x)-(\nabla f)\circ\phi_z(y)\cdot\nabla\phi_z(y)|}{|x-y|^{q-1}}\\
&\leq\sup_{x\not= y}\frac{|(\nabla f)\circ\phi_z(x)-(\nabla f)\circ\phi_z(y)|\cdot \|\nabla\phi_z\|_\infty}{|x-y|^{q-1}}\\
&\qquad+\sup_{x\not= y}\frac{\|\nabla f\|_\infty|\nabla\phi_z(x)-\nabla\phi_z(y)|}{|x-y|^{q-1}}\\
&\stackrel{(\ref{JH1})}{\leq} C\|\nabla f\|_{q-1}(1+|z|)\sup_{x\not= y}\frac{|\phi_z(x)-\phi_z(y)|^{q-1}}{|x-y|^{q-1}}\\
&\quad+C\|\nabla f\|_\infty(1+|z|)^{2-q}\sup_{x\not= y}\frac{|\nabla\phi_z(x)-\nabla\phi_z(y)|^{q-1}}{|x-y|^{q-1}}\\
&\stackrel{(\ref{JH1})}{\leq} C[\nabla f]_{q-1}(1+|z|)^{q}+C\|\nabla f\|_\infty(1+|z|).
\end{align*}
Thus,
$$
[\cL_1 f]_{q-1}\leq C\|f\|_{q-1,\infty}\int_{|z|\geq 1}(1+|z|^{q})\nu(\dif z).
$$
For $q\geq 2$, it follows by similar calculations.
\end{proof}
Let $\cT_t$ be the semigroup associated with $\cL^\nu_{\sigma,b}$.
For any $\varphi\in C^\infty_b(\mR^d)$, as above by Duhamel's formula, we have
\begin{align}
\cT_t\varphi(x)=\cT^0_t\varphi(x)+\int^t_0\cT^0_{t-s}\cL_1\cT_s\varphi(x)\dif s.\label{NJ5}
\end{align}
\bl\label{Le55}
Let $\gamma_{100}$ be as in Theorem \ref{Th45}. If $\int_{|z|\geq 1}|z|^q\nu(\dif z)<\infty$ for some $q>0$, then for any $\beta\in (0,\frac{1}{\gamma_{100}}\wedge 1)$ with 
$q+\beta$ being not an integer, there exists a constant $C>0$ such that for all $t\in(0,1)$ and $\varphi\in L^\infty(\mR^d)$,
\begin{align*}
\|\cT_t\varphi\|_{q+\beta,\infty}\leq C t^{-(q+\beta)\gamma_{100}}\|\varphi\|_\infty.
\end{align*}
\el
\begin{proof}
Fix an irrational number $q_0\in(0,q]$ and choose $m\in\mN$ being large so that 
$$
\theta:=\tfrac{q_0}{m}<\tfrac{1}{\gamma_{100}}\wedge 1.
$$ 
By (\ref{NJ5}), (\ref{ENL}) and Lemma \ref{Le53}, we have
\begin{align*}
\|\cT_t\varphi\|_{m\theta,\infty}&\leq\|\cT^0_t\varphi\|_{m\theta,\infty}+\int^t_0\|\cT^0_{t-s}\cL_1\cT_s\varphi\|_{m\theta,\infty}\dif s\\
&\leq Ct^{-\theta\gamma_{100}}\|\varphi\|_{(m-1)\theta,\infty}+C\int^t_0\|\cT_s\varphi\|_{m\theta,\infty}\dif s,
\end{align*}
which, by Gronwall's inequality, yields that for all $t\in(0,1)$,
$$
\|\cT_t\varphi\|_{m\theta,\infty}\leq Ct^{-\theta\gamma_{100}}\|\varphi\|_{(m-1)\theta,\infty}.
$$
Since $j\theta$ is not an integer for any $j\in\mN$, by iteration we obtain
$$
\|\cT_t\varphi\|_{q_0,\infty}=\|\cT_t\varphi\|_{m\theta,\infty}\leq Ct^{-m\theta\gamma_{100}}\|\varphi\|_{\infty}=Ct^{-q_0\gamma_{100}}\|\varphi\|_{\infty}.
$$
Next we choose $\theta_0\in(0,\frac{1}{\gamma_{100}}\wedge 1)$ and an irrational number $q_0\leq q$ 
so that $q_0+\theta_0=q+\beta$. As above, we have
\begin{align*}
\|\cT_t\varphi\|_{q_0+\theta_0,\infty}
&\leq Ct^{-\theta_0\gamma_{100}}\|\varphi\|_{q_0,\infty}+C\int^t_0(t-s)^{-\theta_0\gamma_{100}}\|\cL_1\cT_s\varphi\|_{q_0,\infty}\dif s\\
&\leq Ct^{-\theta_0\gamma_{100}}\|\varphi\|_{q_0,\infty}+C\|\varphi\|_{q_0,\infty}\int^t_0(t-s)^{-\theta_0\gamma_{100}}\dif s\\
&\leq C(t^{-\theta_0\gamma_{100}}+t^{1-\theta_0\gamma_{100}})\|\varphi\|_{q_0,\infty}.
\end{align*}
Thus,
$$
\|\cT_{2t}\varphi\|_{q_0+\theta_0,\infty}\leq Ct^{-\theta_0\gamma_{100}}\|\cT_t\varphi\|_{q_0,\infty}\leq Ct^{-(q_0+\theta_0)\gamma_{100}}\|\varphi\|_{\infty}.
$$
The proof is complete.
\end{proof}
\bl\label{Le51}
Let $\gamma_{010}$ be as in Theorem \ref{Th45}.
For any $\theta\in(0,1/\gamma_{010}\wedge 1)$, there exists a constant $C>0$ such that for all $t\in(0,1)$ and $\varphi\in C^\infty_b(\mR^d)$,
$$
\|\cT_t\Delta^{\frac{\theta}{2}}\varphi\|_\infty\leq Ct^{-\theta\gamma_{010}}\|\varphi\|_\infty.
$$
\el
\begin{proof}
First of all, we show that
\begin{align}
\|\cT^0_t\Delta^{\frac{\theta}{2}}\varphi\|_\infty\leq C t^{-\theta\gamma_{010}}\|\varphi\|_\infty.\label{NJ1}
\end{align}
Notice that
\begin{align}
\cT^0_t\Delta^{\frac{\theta}{2}}\varphi(x)=\mE\int_{\mR^d}\frac{\varphi(X_t(x)+z)-\varphi(X_t(x))}{|z|^{d+\theta}}\dif z=I_1(x)+I_2(x),\label{NJ2}
\end{align}
where
\begin{align*}
I_1(x)&:=\mE\int_{|z|\leq t^{\gamma_{010}}}\frac{\varphi(X_t(x)+z)-\varphi(X_t(x))}{|z|^{d+\theta}}\dif z,\\
I_2(x)&:=\mE\int_{|z|>t^{\gamma_{010}}}\frac{\varphi(X_t(x)+z)-\varphi(X_t(x))}{|z|^{d+\theta}}\dif z.
\end{align*}
For $I_1(x)$, setting $\varphi_{sz}(x):=\varphi(x+sz)$, we have
\begin{align*}
I_1(x)&=\mE\int_{|z|\leq t^{\gamma_{010}}}\left(\int^1_0z\cdot\nabla\varphi(X_t(x)+sz)\dif s\right)\frac{\dif z}{|z|^{d+\theta}}\\
&=\int_{|z|\leq t^{\gamma_{010}}}\left(\int^1_0z\cdot\cT^0_t\nabla\varphi_{sz}(x)\dif s\right)\frac{\dif z}{|z|^{d+\theta}}.
\end{align*}
Hence,
\begin{align}
\|I_1\|_\infty&\leq\int_{|z|\leq t^{\gamma_{010}}}\left(\int^1_0\|\cT^0_t\nabla\varphi_{sz}\|_\infty\dif s\right)\frac{|z|\dif z}{|z|^{d+\theta}}\no\\
&\stackrel{(\ref{NJ6})}{\leq} Ct^{-\gamma_{010}}\|\varphi\|_\infty\int_{|z|\leq t^{\gamma_{010}}}\frac{\dif z}{|z|^{d+\theta-1}}\leq C t^{-\theta\gamma_{010}}\|\varphi\|_\infty.\label{NJ3}
\end{align}
For $I_2(x)$, we have
\begin{align}
\|I_2\|_\infty\stackrel{(\ref{IU7})}{\leq} C\|\varphi\|_\infty\int_{|z|>t^{\gamma_{010}}}\frac{1}{|z|^{d+\theta}}\dif z\leq C t^{-\theta\gamma_{010}}\|\varphi\|_\infty.\label{NJ4}
\end{align}
Combining (\ref{NJ2})-(\ref{NJ4}), we obtain (\ref{NJ1}). Now, by (\ref{NJ5}) and (\ref{NJ1}), we have
\begin{align*}
\|\cT_t\Delta^{\frac{\theta}{2}}\varphi\|_\infty&\leq\|\cT^0_t\Delta^{\frac{\theta}{2}}\varphi\|_\infty
+\int^t_0\|\cT^0_{t-s}\cL_1\cT_s\Delta^{\frac{\theta}{2}}\varphi\|_\infty\dif s\\
&\leq Ct^{-\theta\gamma_{010}}\|\varphi\|_\infty+\int^t_0\|\cT_s\Delta^{\frac{\theta}{2}}\varphi\|_\infty\dif s,
\end{align*}
which in turn gives the desired estimate by Gronwall's inequality.
\end{proof}
\begin{proof}[Proof of Theorem \ref{Main1}]
Let $X_t(x)$ solve SDE (\ref{SDE1}). By Lemma \ref{Le51} and Lemma \ref{Le64} in Appendix,  
there exists a function $\rho_t(x,y)\in (L^1\cap L^p)(\mR^d)$ for some $p>1$ such that
for all $\varphi\in C_0(\mR^d)$,
$$
\cT_t\varphi(x)=\mE\varphi(X_t(x))=\int_{\mR^d}\varphi(y)\rho_t(x,y)\dif y.
$$
By a further approximation, the above equality also holds for any $\varphi\in L^\infty(\mR^d)$. 
The $q+\eps$-order H\"older continuity of $x\mapsto\cT_t\varphi(x)$
follows by Lemma \ref{Le55}.
\end{proof}
\section{Appendix}

\subsection{Proof of Lemma \ref{Le41}}
Let $\delta$ be as in (\ref{Con}). For $0<\eps<\delta$, let $X^\eps_t(x)=X^\eps_t$ solve the following SDE:
\begin{align}
X^\eps_t=x+\int^t_0b(X^\eps_s)\dif s+\int^t_0\!\!\!\int_{\Gamma^\delta_\eps}\sigma(X^\eps_{s-},z)\tilde N(\dif s,\dif z),\label{SDE99}
\end{align}
where $\Gamma^\delta_\eps:=\{z\in\mR^d: \eps\leq |z|<\delta\}$.
We first prove the following limit theorem.
\bl\label{Le71}
Under {\bf (H$^\sigma_b$)}, there exist a subsequence $\eps_k\to 0$ and a null set $\Omega_0$ such that for all $\omega\notin\Omega_0$,
$$
\lim_{k\to\infty}\sup_{|x|\leq R}\sup_{t\in[0,1]}|X^{\eps_k}_t(x,\omega)-X_t(x,\omega)|=0,\ \ \forall R\in\mN.
$$
\el
\begin{proof}
Set $Z^\eps_t:=X^\eps_t-X_t$. By Burkholder's inequality (\ref{BT2}) and {\bf (H$^\sigma_b$)}, we have for any $p\geq 2$,
\begin{align*}
\mE\left(\sup_{s\in[0,t]}|Z^\eps_s|^p\right)&\leq C\mE\left(\int^t_0|b(X^\eps_s)-b(X_s)|\dif s\right)^p
+C\mE\left(\sup_{t'\in[0,t]}\left|\int^{t'}_0\!\!\!\int_{\Gamma^\eps_0}\sigma(X_{s-},z)\tilde N(\dif s,\dif z)\right|^p\right)\\
&\quad+C\mE\left(\sup_{t'\in[0,t]}\left|\int^{t'}_0\!\!\!\int_{\Gamma^\delta_\eps}(\sigma(X^\eps_{s-},z)-\sigma(X_{s-},z))\tilde N(\dif s,\dif z)\right|^p\right)\\
&\leq C\int^t_0\mE|Z^\eps_s|^p\dif s+C\int_{\Gamma^\eps_0}|z|^p\nu(\dif z)+C\left(\int_{\Gamma^\eps_0}|z|^2\nu(\dif z)\right)^{\frac{p}{2}},
\end{align*}
where $C$ is independent of $\eps,t\in(0,1)$ and $x\in\mR^d$. By Gronwall's inequality, we obtain
$$
\lim_{\eps\to 0}\sup_{x\in\mR^d}\mE\left(\sup_{t\in[0,1]}|X^\eps_t(x)-X_t(x)|^p\right)=0.
$$
Similarly, we can prove that for any $p\geq 2$,
$$
\lim_{\eps\to 0}\sup_{x\in\mR^d}\mE\left(\sup_{t\in[0,1]}|\nabla X^\eps_t(x)-\nabla X_t(x)|^p\right)=0.
$$
Thus, for any $R>0$, by Sobolev's embedding theorem, we have
\begin{align*}
\lim_{\eps\to 0}\mE\left(\sup_{|x|<R}\sup_{t\in[0,1]}|X^\eps_t(x)-X_t(x)|^p\right)
\leq C\lim_{\eps\to 0}\mE\left(\sup_{t\in[0,1]}\|X^\eps_t(\cdot)-X_t(\cdot)\|^p_{\mW^{1,p}(B_R)}\right)=0,
\end{align*}
where $p>d$ and $\mW^{1,p}(B_R)$ is the first order Sobolev space over $B_R:=\{x\in\mR^d: |x|<R\}$. The desired limit follows by a suitable choice of subsequence $\eps_k$.
\end{proof}

Define $\phi(x,z):=x+\sigma(x,z)$. By (\ref{Con}), the mapping $x\mapsto \phi(x,z)$ is invertible for each $|z|\leq\delta$.
Let $\phi^{-1}(x,z)$ be the inverse of $x\mapsto\phi(x,z)$. Write
\begin{align}
\hat\sigma(x,z):=\sigma(\phi^{-1}(x,z),z),\ \ |z|\leq\delta\label{IU1}
\end{align}
and
\begin{align}
\hat b(x)&:=b(x)+\int_{\Gamma^\delta_0}[\sigma(\phi^{-1}(x,z),z)-\sigma(x,z)]\nu(\dif z),\label{IU2}\\
\hat b_\eps(x)&:=b(x)+\int_{\Gamma^\delta_\eps}[\sigma(\phi^{-1}(x,z),z)-\sigma(x,z)]\nu(\dif z).\label{IU3}
\end{align}
By the chain rule, the following lemma is easy.
\bl\label{Le72}
Under {\bf (H$^\sigma_b$)}, there exists a constant $C>0$ such that for all $x\in\mR^d$ and $|z|\leq\delta$,
$$
|\hat\sigma(x,z)|\leq C|z|,\ \ |\nabla_x\hat\sigma(x,z)|\leq C|z|.
$$
Moreover, $\hat b, \hat b_\eps\in C^1_b(\mR^d)$ and for some $C>0$,
$$
\|\hat b_\eps-\hat b\|_\infty+\|\nabla\hat b_\eps-\nabla\hat b\|_\infty\leq C\int_{\Gamma^\eps_0}|z|^2\nu(\dif z).
$$
\el
Fix $T\in[0,1]$. For $t\in[0,T]$, define
$$
\hat L^T_t:=L_{T-}-L_{T-t+}\mbox{ with }L_{T-t+}:=\lim_{s\downarrow t}L_{T-s}.
$$
In particular, $(\hat L^T_t)_{t\in[0,T]}$ is still a L\'evy process with the same L\'evy measure $\nu$ and
\begin{align}
\Delta \hat L^T_t=\Delta L_{T-t}.\label{HG1}
\end{align}
Let ${\hat N}^T(\dif s,\dif z)$ be the Poisson random measure associated with $\hat L^T_t$, i.e.,
$$
\hat N^T((0,t]\times E):=\sum_{0<s\leq t}1_E(\Delta \hat L^T_s),\ \ E\in\sB(\mR^d),
$$
and $\tilde {\hat N}^T(\dif s,\dif z):=\hat N^T(\dif s,\dif z)-\dif s\nu(\dif z)$ the compensated Poisson random measure.
We have
\bl
Let $\hat X^T_t(x)=\hat X^T_t$ solve the following SDE:
\begin{align}
\hat X^T_t=x-\int^t_0\hat b(\hat X^T_s)\dif s-\int^t_0\!\!\!\int_{\Gamma^\delta_0}\hat\sigma(\hat X^T_{s-},z)\tilde {\hat N}^T(\dif s,\dif z),\label{SDE91}
\end{align}
where $\hat\sigma$ and $\hat b$ are defined by (\ref{IU1}) and (\ref{IU2}) respectively. Then 
\begin{align}
\hat X^T_T(x)=X^{-1}_T(x),\ \forall x\in\mR^d,\ a.s.\label{IU5}
\end{align}
\el
\begin{proof} For $\eps\in(0,\delta)$, let $X^\eps_t(x)=X^\eps_t$ solve the following random ODE:
\begin{align*}
X^\eps_t&=x+\int^t_0b(X^\eps_s)\dif s+\int^t_0\!\!\!\int_{\Gamma^\delta_\eps}\sigma(X^\eps_{s-},z)\tilde N(\dif s,\dif z)\\
&=x+\int^t_0\tilde b_\eps(X^\eps_s)\dif s+\sum_{0<s\leq t}\sigma(X^\eps_{s-},\Delta L_s)1_{\Gamma^\delta_\eps}(\Delta L_s),
\end{align*}
where
$$
\tilde b_\eps(x)=b(x)-\int_{\Gamma^\delta_\eps}\sigma(x,z)\nu(\dif z).
$$
By the change of variables, we have
\begin{align*}
X^\eps_{T-t}&=X^\eps_T-\int^T_{T-t}\tilde b_\eps(X^\eps_s)\dif s-\sum_{T-t<s\leq T}\sigma(X^\eps_{s-},\Delta L_s)1_{\Gamma^\delta_\eps}(\Delta L_s)\\
&=X^\eps_T-\int^t_0\tilde b_\eps(X^\eps_{T-s})\dif s-\sum_{0\leq s<t}\sigma(X^\eps_{(T-s)-},\Delta L_{T-s})1_{\Gamma^\delta_\eps}(\Delta L_{T-s}).
\end{align*}
Noticing that if $\Delta L_t\in \Gamma^\delta_\eps$, then 
$$
X^\eps_t-X^\eps_{t-}=\sigma(X^\eps_{t-},\Delta L_t)\Rightarrow X_{t-}=\phi^{-1}(X^\eps_t,\Delta L_t),
$$
and since $\Delta L_T=0$ almost surely, we further have
\begin{align*}
X^\eps_{T-t}&=X^\eps_T-\int^t_0\tilde b_\eps(X^\eps_{T-s})\dif s-\sum_{0<s<t}\hat\sigma(X^\eps_{T-s},\Delta  L_{T-s})1_{\Gamma^\delta_\eps}(\Delta L_{T-s})\\
&\stackrel{(\ref{HG1})}{=}X^\eps_T-\int^t_0\tilde b_\eps(X^\eps_{T-s})\dif s-\sum_{0<s<t}\hat\sigma(X^\eps_{T-s},\Delta \hat L^T_s)1_{\Gamma^\delta_\eps}(\Delta \hat L^T_s),
\end{align*}
where $\hat \sigma(x,z)$ is defined by (\ref{IU1}). In particular,
\begin{align*}
X^\eps_{T-t+}&=X^\eps_T-\int^t_0\tilde b_\eps(X^\eps_{T-s})\dif s-\int^t_0\!\!\!\int_{\Gamma^\delta_\eps}\hat\sigma(X^\eps_{T-s},z){\hat N}(\dif s,\dif z)\\
&=X^\eps_T-\int^t_0\hat b_\eps(X^\eps_{T-s})\dif s-\int^t_0\!\!\!\int_{\Gamma^\delta_\eps}\hat\sigma(X^\eps_{T-s},z)\tilde{\hat N}(\dif s,\dif z),
\end{align*}
where $\hat b_\eps(x)$ is defined by (\ref{IU3}).
On the other hand, let $\hat X^{T,\eps}_t(x)=\hat X^{T,\eps}_t$ solve the following SDE:
$$
\hat X^{T,\eps}_t=x-\int^t_0\hat b_\eps(\hat X^{T,\eps}_s)\dif s-\int^t_0\!\!\!\int_{\Gamma^\delta_\eps}\hat\sigma(\hat X^{T,\eps}_{s-},z)\tilde {\hat N}^T(\dif s,\dif z).
$$
By the uniqueness of solutions to random ODEs, we have
$$
X^\eps_{T-t+}(x)=\hat X^{T,\eps}_t(X^\eps_T(x)),\ \forall x\in\mR^d,\ a.s.
$$
In particular, 
\begin{align}
x=\hat X^{T,\eps}_T(X^\eps_T(x)),\ \forall x\in\mR^d,\ a.s.\label{IU4}
\end{align}
By Lemmas \ref{Le72}, \ref{Le71} and  taking limits for (\ref{IU4}), we obtain
$$
x=\hat X^T_T(X_T(x)),\ \forall x\in\mR^d,\ a.s.
$$
The proof is complete.
\end{proof}
Now we can give
\begin{proof}[Proof of Lemma \ref{Le41}]
By equation (\ref{SDE91}) and a standard calculation, we have for any $p\geq 2$,
$$
\sup_{T\in[0,1]}\sup_{x\in\mR^d}\mE|\nabla \hat X^T_T(x)|^p<\infty,
$$
which, together with (\ref{IU5}), implies that
$$
\sup_{T\in[0,1]}\sup_{x\in\mR^d}\mE(\det(\nabla X^{-1}_T(x)))<\infty.
$$
The desired estimate (\ref{IU7}) then follows by the change of variables and the above estimate.
\end{proof}
\subsection{A criterion for the existence of density}
\bl\label{Le64}
Let $\cT$ be a bounded linear operator in $C_b(\mR^d)$. Assume that for some $\theta\in(0,1)$ and any $\varphi\in C^\infty_b(\mR^d)$,
\begin{align}
\|\cT\Delta^{\frac{\theta}{2}}\varphi\|_\infty\leq C_\theta\|\varphi\|_\infty.\label{64}
\end{align}
Then there exists a measurable function $\rho(x,y)$ with $\rho(x,\cdot)\in (L^1\cap L^p)(\mR^d)$ for some $p>1$ and such that
for any $\varphi\in C_0(\mR^d)$,
\begin{align}
\cT \varphi(x)=\int_{\mR^d}\varphi(y)\rho(x,y)\dif y.\label{HG4}
\end{align}
\el
\begin{proof}
By Riesz's representation theorem, there exists a family of finite signed measures $\mu_x(\dif y)$ such that $x\mapsto\mu_x(\dif y)$ is weakly continuous and
for any $\varphi\in C_0(\mR^d)$,
\begin{align}
\cT\varphi(x)=\int_{\mR^d}\varphi(y)\mu_x(\dif y).\label{JG4}
\end{align}
Let $\varrho$ be a nonnegative symmetric smooth function with compact support and $\int_{\mR^d}\varrho(y)\dif y=1$. Let $\varrho_\eps(y):=\eps^{-d}\varrho(\eps^{-1}y)$ 
be a family of mollifies. For $R>0$, let $\chi_R:\mR^d\to[0,1]$ be a smooth cutoff function with 
$$
\chi_R(x)=1,\ \ |x|\leq R,\ \ \chi_R(x)=0,\ \ |x|\geq 2R.
$$
For $\varphi\in L^\infty(\mR^d)$,  set 
$$
\varphi_\delta(x):=\varphi*\varrho_\delta(x), \ \ \varphi^R_{\delta,\eps}(x):=(\varphi_\delta\chi_R)*\varrho_\eps(x)
$$
and
$$
\mu^\eps_x(z):=\int_{\mR^d}\varrho_\eps(y-z)\mu_x(\dif y).
$$
It is easy to see that $\mu^\eps_x\in \cap_{k}\mW^{k,1}(\mR^d)$ and $\Delta^{\frac{\theta}{2}}\varphi^R_{\delta,\eps}\in C_0(\mR^d)$. 
Thus, by (\ref{JG4}) we have
\begin{align*}
\cT\Delta^{\frac{\theta}{2}}\varphi^R_{\delta,\eps}(x)&=\int_{\mR^d}(\Delta^{\frac{\theta}{2}}(\varphi_\delta\chi_R))*\varrho_\eps(y)\mu_x(\dif y)\\
&=\int_{\mR^d}\Delta^{\frac{\theta}{2}}(\varphi_\delta\chi_R)(z)\mu^\eps_x(z)\dif z=\int_{\mR^d}\varphi_\delta\chi_R(z)\Delta^{\frac{\theta}{2}}\mu^\eps_x(z)\dif z,
\end{align*}
which yields by (\ref{64}) that
$$
\left|\int_{\mR^d}\varphi_\delta\chi_R(z)\Delta^{\frac{\theta}{2}}\mu^\eps_x(z)\dif z\right|\leq C_\theta\|\varphi^R_{\delta,\eps}\|\leq C_\theta\|\varphi\|_\infty.
$$
Letting $R\to\infty$ and $\delta\to 0$, by the dominated convergence theorem, we obtain that for all $\varphi\in L^\infty(\mR^d)$,
$$
\left|\int_{\mR^d}\varphi(z)\Delta^{\frac{\theta}{2}}\mu^\eps_x(z)\dif z\right|\leq C_\theta\|\varphi\|_\infty,
$$
which gives
$$
\sup_{x\in\mR^d}\sup_{\eps\in(0,1)}\|\Delta^{\frac{\theta}{2}}\mu^\eps_x\|_1\leq C_\theta.
$$
Moreover, we also have
$$
\sup_{x\in\mR^d}\sup_{\eps\in(0,1)}\|\mu^\eps_x\|_1\leq \sup_{\|\varphi\|_\infty\leq 1}\|\cT\varphi\|_\infty.
$$
By Sobolev's embedding theorem, there is a $p>1$ such that
$$
\sup_{x\in\mR^d}\sup_{\eps\in(0,1)}\|\mu^\eps_x\|_p<\infty.
$$
Since $L^p(\mR^d)$ is reflexive, for each fixed $x\in\mR^d$, there is a subsequence $\eps_k\to 0$ and
$\rho(x,\cdot)\in (L^1\cap L^p)(\mR^d)$ such that for any $\varphi\in C_0(\mR^d)\subset L^\frac{p}{p-1}(\mR^d)$,
\begin{align}
\cT\varphi_{\eps_k}(x)=\int_{\mR^d}\mu^{\eps_k}_x(z)\varphi(z)\dif z\stackrel{k\to\infty}{\to}\int_{\mR^d}\rho(x,z)\varphi(z)\dif z.\label{HG5}
\end{align}
On the other hand, for any $\varphi\in C_0(\mR^d)$,
$$
\|\cT\varphi_\eps-\cT\varphi\|_\infty\leq C\|\varphi_\eps-\varphi\|_\infty\stackrel{\eps\to 0}{\to} 0,
$$
which together with (\ref{HG5}) yields (\ref{HG4}).
\end{proof}

{\bf Acknowledgements:}

The author is grateful to Linlin Wang and Longjie Xie for their useful conversations. This work is supported by NNSFs of China (Nos. 11271294, 11325105).

\end{document}